\documentclass[twoside, 12pt]{article}

\usepackage{amsmath, amsthm, amssymb, amscd, mathptmx}

\usepackage{amsfonts}

\theoremstyle{plain}
\newtheorem{Thm}{\sc Theorem}[section]
\newtheorem{Theorem}[Thm]{\sc Theorem}
\newtheorem{Corollary}[Thm]{\sc Corollary}

\newtheorem*{Corollary*}{\sc Corollary}

\newtheorem{Proposition}[Thm]{\sc Proposition}
\newtheorem*{Proposition*}{\sc Proposition}
\newtheorem{Lemma}[Thm]{\sc Lemma}

\theoremstyle{definition}
\newtheorem{Definition}[Thm]{Definition}

\theoremstyle{remark}
\newtheorem{Remark}[Thm]{Remark}
\newtheorem{Example}[Thm]{Example}
\newtheorem*{Example*}{Example}
\newtheorem*{Remark*}{Remark}

\newcommand{\CC}{{\mathbb C}}

\newcommand{\NN}{{\mathbb N}}
\newcommand{\ZZ}{{\mathbb Z}}

\newcommand{\PP}{{\mathbb P}}
\newcommand{\QQ}{{\mathbb Q}}
\newcommand{\RR}{{\mathbb R}}

\newcommand{\cA}{{\mathcal A}}

\newcommand{\cB}{{\mathcal B}}

\newcommand{\cS}{{\mathcal S}}
\newcommand{\cE}{{\mathcal E}}
\newcommand{\cF}{{\mathcal F}}
\newcommand{\cG}{{\mathcal G}}
\newcommand{\cH}{{\mathcal H}}

\newcommand{\cK}{{\mathcal K}}
\newcommand{\cL}{{\mathcal L}}
\newcommand{\cM}{{\mathcal M}}

\newcommand{\cO}{{\mathcal O}}

\newcommand{\cV}{{\mathcal V}}

\newcommand{\Pic}{{\mathop{\rm Pic \, }}}

\newcommand{\Refl}{{\mathop{\operatorname{Ref}\, }}}

\newcommand{\ch}{{\mathop{\rm ch \, }}}

\newcommand{\id}{\mathop{\rm Id}}

\newcommand{\im}{\mathop{\rm im \, }}

\newcommand{\rk}{{\mathop{\rm rk \,}}}

\newcommand{\Sym}{{\mathop{{\rm Sym \, }}}}

\newcommand{\Spec}{{\mathop{{\rm Spec\, }}}}

\newcommand{\Vect}{{\mathop{{\rm Vect \,}}}}

\newcommand{\Sing}{{\mathop{{\rm Sing \,}}}}

\newcommand{\refl}{{\mathop{\operatorname{\rm ref} }}}

\usepackage[T1]{fontenc}
\usepackage[hidelinks]{hyperref}
\newcommand{\cQ}{{\mathcal Q}}

\allowdisplaybreaks[1]

\begin{document}
\markboth{\rm }{}

\title{Intersection theory and Chern classes on normal varieties II}
\author{Adrian Langer}
\date{November 16, 2026}

\maketitle

{\noindent \sc Address:}\\
Institute of Mathematics, University of Warsaw,
ul.\ Banacha 2, 02-097 Warszawa, Poland\\
e-mail: {\tt alan@mimuw.edu.pl}

\medskip

\begin{abstract}
	We prove a uniform bound linear in the rank for the local
	correction terms in the Riemann--Roch formula for reflexive
	sheaves on normal surfaces over algebraically closed fields
	of arbitrary characteristic. We also establish the usual tensor
	and symmetric-power formulas for second Chern characters
	and prove Wahl's asymptotic Riemann--Roch conjecture and its
	higher-rank generalization. These results yield second Chern
	classes of reflexive sheaves on arbitrary normal projective
	varieties, compatible with pullback by finite surjective
	morphisms. We characterize their discriminants as infima of
	normalized discriminants of compatible vector bundles on
	smooth alterations.  As applications we prove Bogomolov's inequality, effective restriction and boundedness theorems on normal projective varieties, a Miyaoka-type inequality,
	 and a Segre-class criterion for
	bigness of reflexive sheaves on normal surfaces.
\end{abstract}

\section*{Introduction}

In \cite{Wa} Wahl introduced relative Chern classes of rank $2$ vector bundles on resolutions of complex surface singularities. In \cite{La0} and \cite{LaI} we extended these classes to arbitrary rank and characteristic, and used them to construct Chern classes of reflexive sheaves on normal surfaces over an algebraically closed field. The most important feature of these classes is that they behave like orbifold Chern classes and they are useful in asymptotic Riemann--Roch computations. For a normal proper surface $X$ and a reflexive
coherent $\cO_X$-module $\cE$ of rank $r$, these classes satisfy also the Riemann--Roch formula
$$\chi(X,\cE)=\frac12c_1(\cE).(c_1(\cE)-K_X)
-\int_Xc_2(\cE)+r\chi(X,\cO_X)
+\sum_{x\in\Sing X}a(x,\cE),$$
where the terms $a(x,\cE)$ depend only on the local reflexive modules at
the singular points (see \cite[Theorem 4.4]{LaI}). The main aim of this paper is to prove that
for every closed point $x\in X$ there exists a constant $C_x$ such that $|a(x,\cE)|\le C_xr$ for all $\cE$.
This result was conjectured in \cite[Conjecture 8.2]{La0} (see also \cite[Remark  3.14]{LaI}).
It implies Wahl's main conjecture \cite[Introduction]{Wa} on the local asymptotic Riemann--Roch formula  
(see also \cite[Conjecture 4.1]{La0} and its generalization \cite[Conjecture 3.18]{LaI}). As a corollary we obtain the following theorem improving \cite[Theorem 0.3]{LaI}.

\begin{Theorem}\label{main1} 
	Let $X$ be a normal proper algebraic surface defined over an algebraically closed field.
	Then there exists a constant $C$ depending only on $X$ such that for every rank $r$  coherent reflexive $\cO_X$-module $\cE$ on $X$ we have
	$$\left|\chi (X, \cE) -\left(\frac {1}{2}c_1 (\cE). (c_1(\cE) -K_X)- \int_Xc_2 (\cE) +r\chi (X, \cO_X)\right)\right|\le C r.$$
\end{Theorem}

We also prove that the Chern character of reflexive sheaves on a normal projective surface is a homomorphism of rings, as conjectured in \cite[Section 5.5]{LaI} (see Theorem~\ref{surface:chern-ring}).
More precisely, let $K^{\refl}(X)$ be the Grothendieck ring of reflexive
sheaves, with relations coming from locally split short exact sequences
and multiplication given by reflexive tensor product. Then
$$\ch:K^{\refl}(X)\longrightarrow A^*(X)\otimes\RR,
\qquad [\cE]\longmapsto\rk\cE+c_1(\cE)+\ch_2(\cE)$$
is a ring homomorphism, where the target has the ring structure induced
by Mumford's intersection product.

We use Theorem \ref{main1} to construct a good theory for the second Chern class (or character) for reflexive sheaves in higher dimensions. This was already done in \cite{LaI} in  positive characteristic or for varieties with at most quotient singularities in codimension two in characteristic zero. Here we extend these results to arbitrary  normal projective varieties and we prove the following theorem strengthening  \cite[Theorem 0.4]{LaI} (see Theorems \ref{high:properties} and \ref{thm:c_2-mult-by-degree}).

\begin{Theorem}
For any normal projective variety $X$ of dimension $n\ge 2$ over an algebraically closed field $k$ and for any 
coherent reflexive $\cO_X$-module $\cE$ on $X$ there exists a $\ZZ$-multilinear symmetric 
form $\int _X c _2 (\cE): \Pic X ^{\times (n-2)}\to \RR$ such that:
	\begin{enumerate}
		\item If $\cE$ is a vector bundle on $X$ then $$\int _X c _2 (\cE)L_1...L_{n-2}=\int_X c _2 (\cE)\cap c_1(L_1)\cap ...\cap c_1(L_{n-2})\cap[X].$$
		\item If $k\subset K$ is an algebraically closed field extension then 
		$$\int _{X_K} c _2 (\cE_K)(L_1)_K...(L_{n-2})_K=\int _X c _2 (\cE)L_1...L_{n-2}.$$
		\item If $n>2$ and  $ L_1$ is very ample then for a very general hypersurface $H\in |L_1|$ we have
		$$\int _X c_2 (\cE)L_1...L_{n-2}=\int _{H}c_2 (\cE|_{H})L_2|_{H}...L_{n-2}|_{H}.$$
	\item If $\pi:Y\to X$ is a finite surjective morphism from a normal
	projective variety $Y$, then
	$$\int_Yc_2(\pi^{[*]}\cE)\pi^*L_1...\pi^*L_{n-2}
	=\deg\pi\cdot\int_Xc_2(\cE)L_1...L_{n-2}.$$
	\end{enumerate}
\end{Theorem}

The second Chern class satisfies also the following asymptotic version of the Riemann--Roch theorem (see Proposition~\ref{high:limit} and Definition~\ref{high:definition}).

\begin{Theorem} 
Let $\cE$ be a rank $r$ coherent reflexive $\cO_X$-module on $X$ as above. Then for any line bundles $L_1,...,L_{n-2}$ on $X$ we have
	\begin{align*}
		&\chi\left(X,c_1(L_1)...c_1(L_{n-2})\cdot\Sym^{[m]}\cE\right)\\
		&\qquad=\frac{m^{r+1}}{(r+1)!}
		\left(c_1(\cE)^2.L_1...L_{n-2}
		-\int_Xc_2(\cE)L_1...L_{n-2}\right)+O(m^r).
	\end{align*}
\end{Theorem}

We also prove the following characterization of the discriminant $\int_X\Delta(\cE)=2r \int_Xc_2(\cE)-(r-1)c_1(\cE)^2$: 

\begin{Theorem}\label{main4}
	Let $X$ be a normal projective variety of dimension $n\geq2$
	over an algebraically closed field, and let
	$L_1,\ldots,L_{n-2}$ be nef line bundles on $X$.
	For every reflexive sheaf $\cE$ on $X$ we have
	\[
	\int_X\Delta(\cE)L_1\cdots L_{n-2}
	=
	\inf_{f,\cF}
	\frac{
		\int_Z\Delta(\cF)f^*L_1\cdots f^*L_{n-2}
	}{\deg f},
	\]
	where $f:Z\to X$ runs through all alterations from smooth projective varieties, and $\cF$ runs through all
	vector bundles on $Z$ satisfying
	\(
	g_{[*]}\cF\simeq h^{[*]}\cE,
	\)
	with
	\(
	Z\xrightarrow{g}Y\xrightarrow{h}X
	\)
	the Stein factorization of $f$.	
\end{Theorem}

The above theorem easily implies Bogomolov's inequality for (strongly) slope semistable sheaves on normal projective varieties (see Corollary \ref{Bogomolov-inequality}). This inequality implies a strong restriction theorem and boundedness results on normal projective varieties in characteristic zero (see Section \ref{Applications}).  As a corollary we also prove a generalization of Miyaoka's inequality on the second Chern class of the cotangent bundle of a smooth non-uniruled complex variety with nef $K_X$ to normal projective varieties (see Corollary \ref{Miyaoka-cotangent}). 

In the case of normal complex K\"ahler varieties a weaker form of Bogomolov's inequality for slope semistable sheaves was proven by Chen in \cite{Chen2025} (see Subsection \ref{Chen}).

\medskip

The paper is organized as follows. In Section~1 we prove the
uniform linear bound for the local Riemann--Roch correction
terms. In Section~2 we establish the tensor and symmetric-power
formulas for relative Chern characters and derive the local
asymptotic Riemann--Roch formula. Section~3 studies the Chern
character on normal surfaces and constructs second Chern
classes in higher dimensions, proving their compatibility
with finite pullbacks. In Section~4 we characterize
discriminants by approximation with vector bundles on smooth
alterations and compare our invariant with Chen's.
Section~5 contains applications to Bogomolov's inequality,
restriction theorems, boundedness, Miyaoka-type inequalities,
and bigness of reflexive sheaves and finiteness of curves of
low genus on normal surfaces. Finally, the appendix shows
that in higher degrees the expected asymptotic
Riemann--Roch formula is incompatible with the usual duality
formula, already for rank-one reflexive sheaves.

\subsection*{Notation}

We use the notation of \cite{LaI}. For  a coherent sheaf $\cE$ on a normal scheme $X$ 
we write $\Sym^{[m]}\cE=(\Sym^m\cE)^{**}$. For coherent sheaves $\cE$ and $\cF$ 
we also write $\cE\hat\otimes\cF:=(\cE\otimes\cF)^{**}$.
If $f: X\to Y$ is a morphism of normal schemes, $\cE$ is a coherent sheaf on $X$ and $\cF$ a coherent sheaf on $Y$ then we write
$f^{[*]}\cF=(f^*\cF)^{**}$ and $f_{[*]}\cE=(f_*\cE)^{**}$.

\section{Local correction to the Riemann--Roch theorem}

In this section we fix the following notation.
Let $k$ be an algebraically closed field and let $A$ be the
henselization of the local ring of a normal algebraic surface at
a closed point. Put $X=\Spec A$, and write $x$ for its closed
point. We fix a projective good resolution $f:\tilde X\to X$, which is an
isomorphism over $U=X\backslash\{x\}$. We write
$E=\bigcup_{i=1}^sE_i$ for the reduced exceptional divisor. Its
irreducible components are smooth projective curves, and its
intersection matrix is negative definite.

\subsection{Relative Chern classes and the correction term}

Here we recall a few notions from \cite{La0} and \cite{LaI} (see also \cite{Wa}).

For a vector bundle $\cF$ of rank $r$ on $\tilde X$, 
$c_1(f,\cF)$ is the exceptional $\QQ$-divisor determined by
$$c_1(f,\cF).E_i=\deg(\cF|_{E_i}).$$
We use the relative second Chern class $c_2(f,\cF)$ from
\cite[Definition 3.1]{LaI}. More explicitly, let
$\tilde\pi:\tilde Y\to\tilde X$ be a generically finite proper
morphism from a regular surface, and let
$\tilde Y\xrightarrow{g}Y\to X$ be the Stein factorization of
$f\tilde\pi$. For a filtration of $\tilde\pi^*\cF$ with line
bundle quotients $\cL_1,...,\cL_r$, put $L_i=c_1(g,\cL_i)$.
Then
$$c_2(f,\cF)=\inf\frac{\sum_{i<j}L_i.L_j}{\deg\tilde\pi},$$
where the infimum is taken over these morphisms and filtrations.
We also  set
$$\ch_2(f,\cF):=\frac12 c_1(f,\cF)^2-c_2(f,\cF),\quad
\Delta(f,\cF):=2r c_2(f,\cF)-(r-1)c_1(f,\cF)^2.$$
The relative Euler characteristic is defined by
$$\chi(f,\cF)=
\dim H^0(X,f_{[*]}\cF/f_*\cF)
+\dim H^0(X,R^1f_*\cF).$$
We also  define the correction term to a local Riemann--Roch formula as
\begin{equation*}
	a(f,\cF)
	:=\chi(f,\cF)-r\chi(f,\cO_{\tilde X})
	+\frac12  c_1(f,\cF).(c_1(f,\cF)-K_{\tilde X})-c_2(f,\cF).
\end{equation*}
By \cite[Proposition 3.7]{LaI}, this number depends only on
$f_{[*]}\cF$. For a reflexive sheaf $\cE$ on $X$, we denote it by
$a(x,\cE)$. The main aim of this section is to prove the following
bound on the correction terms $a(x,\cE)$.

\begin{Theorem}\label{geom:linear-bound}
	There exists a constant $C_x$ depending only on the germ $(X,x)$
	such that for every coherent reflexive $\cO_X$-module $\cE$
	of rank $r$ we have
	\[
	|a(x,\cE)|\leq C_xr.
	\]
\end{Theorem}

The proof in characteristic zero uses vector bundles minimizing
the relative Euler characteristic. Elementary transformations give
bounds for the Harder--Narasimhan slopes of their restrictions to
the exceptional curves. An induction on rank, using complete
intersections in a projective bundle, converts these bounds into
a relative discriminant estimate. We then construct finite covers
using roots of a fixed collection of functions. The eigensheaves of
these covers belong to a fixed finite collection of exceptional
line bundles. This controls the geometric genera of the covering
singularities and gives the required linear bound.

In positive characteristic the proof is easier and it uses the results of \cite{LaI} and 
a bound on the difference between a Frobenius pullback of a reflexive module and its reflexive hull.

\subsection{Minimal extensions}

\begin{Definition}
	Let $\cE$ be a coherent reflexive $\cO_X$-module.
	We say that a vector bundle $\cF$ on $\tilde X$ is a \emph{minimal
		extension of $\cE$} if $\chi(f,\cF)$ is minimal among the relative Euler
	characteristics of vector bundles with $f_{[*]}\cF\simeq \cE$.
\end{Definition}

Since  $\chi(f,\cF)$ is a non-negative integer, every reflexive sheaf on $X$ admits a minimal extension.

\begin{Lemma}\label{geom:minimal-slopes}
	Let $\cF$ be a minimal extension of some reflexive sheaf. If $C$ is an irreducible component of
	$E$ of genus $g$, then
	\[
	g-1\leq\mu_{\min}(\cF|_C)
	\leq\mu_{\max}(\cF|_C)\leq g-1-C^2.
	\]
\end{Lemma}

\begin{proof}
	Let $\cF|_C\to\cQ$ be the minimal destabilizing quotient and let $\cF'$ be the
	kernel of $\cF\to i_*\cQ$, where $i:C\hookrightarrow\tilde X$. The
	bundle $\cF'$ is locally free and agrees with $\cF$ off $C$. Since their
	reflexive pushforwards are the same, the exact sequence gives
	\[
	\chi(f,\cF')-\chi(f,\cF)=\chi(C,\cQ).
	\]
	Minimality implies that $\chi(C,\cQ)\geq 0$. Hence $\mu_{\min}(\cF|_C)=\mu(\cQ)\ge g-1$.
	
	Similarly, let $\cS\subset\cF|_C$ be the maximal destabilizing subbundle. The elementary
	transformation
	\[
	\cF''=\ker\bigl(\cF(C)\longrightarrow
	i_*((\cF|_C/\cS)\otimes\cO_C(C))\bigr)
	\]
	is locally free and fits into an exact sequence
	\[
	0\longrightarrow\cF\longrightarrow\cF''
	\longrightarrow i_*(\cS(C))\longrightarrow0.
	\]
	Thus
	\[
	\chi(f,\cF'')-\chi(f,\cF)
	=-\chi(C,\cS(C))\geq0,
	\]
	which gives $\mu_{\max}(\cF|_C)=\mu(\cS)\leq g-1-C^2$.
\end{proof}

\subsection{Relative discriminants }

\begin{Lemma}\label{geom:discriminant-spread}
	Assume that the base field $k$ has characteristic zero.
There exists a function $\varphi: \NN\to \RR _{\ge 0}$  such that
for every rank $r$ vector bundle $\cF\in  \Vect (\tilde X)$,  if 
	$T$ is an exceptional $\QQ$-divisor such that for every $i$ we have $	TE_i>0$ and
	\[
	\mu_{\max}(\cF|_{E_i})-\mu_{\min}(\cF|_{E_i})\leq TE_i,
	\]
then
	\(
	\Delta(f,\cF)\leq -\varphi(r) T^2.
	\)
\end{Lemma}

\begin{proof}
	We argue simultaneously for all good resolutions of algebraic normal
surface singularities and construct $\varphi$ recursively, starting with
$\varphi(1)=0$.
	Assume that $r\geq2$. On
	$p:\PP_{\tilde X}(\cF)\to\tilde X$ write
	\[
	\xi=c_1(\cO_{\PP(\cF)}(1)),\qquad
	\eta=\xi-\frac1r p^*c_1(\cF),\qquad
	B=\eta+2p^*T.
	\]
	The slope criterion for ampleness on a curve in characteristic
	zero (see \cite[Section 6.4]{Laz}) shows that
	$B|_{\PP(\cF|_{E_i})}$ is ample: its least slope is
	\[
	\mu_{\min}(\cF|_{E_i})-\mu(\cF|_{E_i})+2TE_i
	\geq TE_i>0.
	\]
	Consequently $B$ is ample relative to $fp$. Choose a sufficiently
large divisible integer $m$ so that $mB$ is an integral very ample class
relative to $X$, and intersect $r-1$ general members of this linear
system. Bertini's theorem on each $\PP(\cF|_{E_i})$, together with
transversality on the fibers over the nodes of $E$, gives a surface
$\tilde Y$ such that each
\[
C_i=\tilde Y\cap p^{-1}(E_i)
\]
is a smooth irreducible curve and the curves $C_i$ have simple
normal crossings. The surface $\tilde Y$ is regular along these
curves. Its nonregular locus is closed and proper over $X$ and
does not meet the closed fiber, so it is empty.

The morphism
\[
\tilde \pi=p|_{\tilde Y}:\tilde Y\longrightarrow\tilde X
\]
is proper and generically finite of degree $d=m^{r-1}$, and each
$C_i$ dominates $E_i$. A positive-dimensional fiber over a point
of $E_i$ would contain a curve component of $C_i$, contradicting
its irreducibility and dominance over $E_i$. Thus $\tilde\pi$ is
quasi-finite above $E$. The locus in $\tilde X$ over which its fibers
have positive dimension is closed, and its image under $f$ is a
closed subset of the local scheme $X$ not containing $x$.
It is therefore empty, so $\tilde\pi$ is finite.

Every irreducible component of the regular surface $\tilde Y$
dominates $\tilde X$. Since these components are disjoint and
each meets $C_i$, irreducibility of $C_i$ implies that $\tilde Y$
is integral. Moreover, $\tilde\pi$ is flat, since both surfaces
are regular. Thus $C_i=\tilde\pi^*E_i$ and $C_i\to E_i$ has
degree $d$. Let
\(
\tilde Y\xrightarrow{g}Y\xrightarrow{\pi} X
\)
be the Stein factorization of $f\tilde\pi$. The exceptional
divisor of $g$ is $\bigcup_i C_i$.

	Restricting the universal quotient gives an exact sequence
	\[
	0\longrightarrow\cK\longrightarrow \tilde \pi^*\cF
	\longrightarrow\cL\longrightarrow0.
	\]
Since	$p_*\eta^{r-1}=1$ and $p_*\eta^r=0$, we have
	\[
	\deg\cL|_{C_i}
	=d\left(\frac{\deg\cF|_{E_i}}r+2(r-1)TE_i\right).
	\]
	This implies that
	\[
	c_1(g,\cL)=\frac1r \tilde \pi^*c_1(f,\cF)+2(r-1)\tilde \pi^*T,
	\]
	and hence
	\begin{equation*}
		c_1(g,\cK)=\frac{r-1}{r}\tilde \pi^*c_1(f,\cF)-2(r-1)\tilde \pi^*T.
	\end{equation*}
	
	Finite pullback of bundles on smooth curves in characteristic zero
	preserves semistability, so $\mu_{\max}(\cK|_{C_i})
	\leq d\mu_{\max}(\cF|_{E_i})$.
	For a vector bundle $V$ of rank $s$ on a smooth curve we have
$s \mu (V)\le 
	(s-1)\mu_{\max}(V)+\mu_{\min}(V)$
and hence
	$$\mu_{\max}(V)-\mu_{\min}(V)
	\leq s\bigl(\mu_{\max}(V)-\mu(V)\bigr).$$ 
Since 	$\mu(\cK|_{C_i})=d\left(\mu(\cF|_{E_i})-2TE_i\right)$, this implies
	\[
	\mu_{\max}(\cK|_{C_i})-\mu_{\min}(\cK|_{C_i})\leq (r-1)(d\mu_{\max}(\cF|_{E_i})-\mu(\cK|_{C_i}))
	\leq3(r-1)dTE_i.
	\]
Applying induction to $\cK$ and $T'=3(r-1)\tilde \pi^*T$ we get
$$\Delta(g,\cK)\le -\varphi(r-1)(T')^2=-9(r-1)^2\varphi(r-1)dT^2.$$
	
	By the definition of the relative second Chern class we have
	\[
	c_2(g,\tilde \pi^*\cF)\leq c_2(g,\cK)+c_1(g,\cK).c_1(g,\cL).
	\]
	This yields
	\begin{align*}
		\Delta(g,\tilde \pi^*\cF)
		&\leq\frac r{r-1}\Delta(g,\cK)
		-\frac{(c_1(g,\cK)-(r-1)c_1(g,\cL))^2}{r-1}\\
		&=\frac r{r-1}\Delta(g,\cK)
		-\frac{(2r(r-1)\tilde \pi^*T)^2}{r-1}\\
		&\le -9r(r-1)\varphi(r-1)dT^2-4r^2(r-1)dT^2.
	\end{align*} 
By \cite[Proposition 3.2, (4)]{LaI} we have $ 	\Delta(g,\tilde \pi^*\cF)=d 	\Delta(f,\cF)$, so 
	\[
	\Delta(f,\cF)
	\leq -\bigl(9r(r-1)\varphi(r-1)+4r^2(r-1)\bigr)T^2.
	\] 
Taking $\varphi (r):=9r(r-1)\varphi(r-1)+4r^2(r-1)$, we get the required assertion.
\end{proof}

\begin{Remark}
The proof shows that  $\varphi (r)$ depends only on $r$ and it is independent of the singularity or its resolution.	
\end{Remark}

\subsection{Finite coverings of surface singularities}

\begin{Proposition}\label{geom:coverings}
	There exist a constant $P_f\geq0$ and, for every sufficiently large prime number $e$ invertible in $k$,
	a finite morphism $\tilde\pi_e:\tilde X_e\to\tilde X$ of degree $e^2$ from a regular
	integral surface, whose composite with $f$ has Stein factorization
	\[
	\tilde X_e\xrightarrow{f_e}X_e\xrightarrow{\pi_e}X,
	\]
	such that the following conditions are satisfied. 
	\begin{enumerate}
		\item $f_e$ is a good resolution. If $C$ runs through the components
		of $\tilde\pi_e^{-1}(E_i)$, then these components are pairwise disjoint
		and
		\(
		\tilde\pi_e^*E_i=e\sum C.
		\)
		\item $\chi (f_e, \cO_{\tilde X_e})\leq P_f\, \deg\pi_e $.
		\item For every vector bundle $\cF$ on $\tilde X$, we have
		\(
		a(f_e,\tilde\pi_e^*\cF)= \deg\pi_e \cdot a(f,\cF).
		\)
	\end{enumerate}
\end{Proposition}

\begin{proof}
	If $f:\tilde X\to X$ is not an isomorphism, choose effective exceptional divisors
$Z_1=\sum_i a_iE_i$ and $ Z_2=\sum_i b_iE_i$ so that all $a_i,b_i$ are positive, the line bundles
	$\cO_{\tilde X}(-Z_1)$ and $\cO_{\tilde X}(-Z_2)$ are $f$-very ample and $a_ib_j-a_jb_i\ne0$ whenever $E_i$ and $E_j$ intersect. Choose general sections of these two line bundles, viewed as functions
	$h_1,h_2$ on $X$, and write
	\[
	\operatorname{div}(h_\nu)=Z_\nu+H_\nu
	\]
for $\nu=1,2.$	
	We can choose them so that the $H_\nu$ are reduced and smooth, have no
	common components, and $E+H_1+H_2$ has simple normal crossings.
	Both $H_\nu$ are nonempty, since $H_\nu.E_i=-Z_\nu.E_i>0$.
	If  $f$ is an isomorphism we set $Z_1=Z_2=0$
	and take $h_1,h_2$ to be a regular system of parameters at $x$, with
	$H_\nu=\operatorname{div}(h_\nu)$.
	
	Let $e$ be a prime invertible in $k$ and dividing none of the nonzero
	integers $a_i$, $b_i$, and $a_ib_j-a_jb_i$ occurring above.
	Let $\tilde X_e$ be the normalization of $\tilde X$ in
	\(
	K(X)(h_1^{1/e},h_2^{1/e}).
	\)
	Valuations along components of $H_1$ and $H_2$ show that the classes of
	$h_1,h_2$ are independent in $K(X)^\times/K(X)^{\times e}$.
	Thus $d_e:=\deg\tilde\pi_e=\deg\pi_e=e^2$,
	and we have the eigensheaf decomposition
	\begin{equation}\label{geom:eigensheaves}
		\tilde\pi_{e*}\cO_{\tilde X_e}
		=\bigoplus_{j,k=0}^{e-1}
		\cO_{\tilde X}(D_{j,k,e}),
		\qquad
		D_{j,k,e}
		=\sum_i\left\lfloor\frac{ja_i+kb_i}e\right\rfloor E_i
	\end{equation}
	(cf.~\cite[Corollary 3.11]{Esnault-Viehweg1992}).
Let us check that $\tilde X_e$ is regular. This can be checked
after passing to strictly henselian local rings. Since $e$ is
invertible, every unit admits an $e$-th root, so unit factors
in $h_1,h_2$ can be omitted when describing the local cover.

At a point lying on only one component of $E+H_1+H_2$, let
$u$ be a local equation of that component. Then $h_1,h_2$
are powers of $u$, up to units, and at least one of their
orders of vanishing is prime to $e$. The local normalization
is therefore a disjoint union of regular covers obtained
by adjoining $u^{1/e}$.

At an intersection of two components, choose regular
parameters $u,v$ defining them. Up to units, the possible
local expressions are
\[
(h_1,h_2)=
\begin{cases}
	(u^{a_i}v^{a_j},u^{b_i}v^{b_j}),
	&\text{at }E_i\cap E_j,\\
	(u^{a_i}v,u^{b_i}),
	&\text{at }E_i\cap H_1,\\
	(u^{a_i},u^{b_i}v),
	&\text{at }E_i\cap H_2,\\
	(u,v),
	&\text{at }H_1\cap H_2.
\end{cases}
\]
In the first case, $a_ib_j-a_jb_i$ is invertible modulo $e$,
so adjoining $e$-th roots of $h_1,h_2$ is equivalent to
adjoining $u^{1/e}$ and $v^{1/e}$. The same conclusion holds
in the second and third cases because $b_i$ and $a_i$,
respectively, are invertible modulo $e$, and it is immediate
in the fourth case. These local covers are regular.
Away from $E+H_1+H_2$, the cover is \'etale.

This proves that $\tilde X_e$ is regular and that the reduced
inverse image of $E$ has simple normal crossings. The same
local descriptions show that the irreducible components of
$\tilde\pi_e^{-1}(E_i)$ are pairwise disjoint and regular, and
that $\tilde\pi_e^*E_i=e\sum C.$
The Stein factorization therefore gives the required good
resolution $f_e$.

	Let $R_e$ be the ramification divisor, so that
	$K_{\tilde X_e}=\tilde\pi_e^*K_{\tilde X}+R_e$.
	By construction, the ramification index is $e$ along every $E_i$ and
	$H_\nu$ (see \cite[Lemma 3.15]{Esnault-Viehweg1992}),
	and we have
	\[
	R_e=\left(1-\frac1e\right) \tilde\pi_e^*\!\left(
	\sum_i E_i
	+H_1+H_2
	\right).
	\]

	All $D_{j,k,e}$ belong to the fixed finite set
	\(
	\mathcal D_f=
	\left\{\sum_i c_iE_i:
	0\leq c_i\leq a_i+b_i,\ c_i\in\ZZ\right\}.
	\)
	So
	\[
	\chi (f_e, \cO_{\tilde X_e})=\sum_{j,k=0}^{e-1}\chi(f,\cO_{\tilde X}(D_{j,k,e}))
	\leq d_eP_f,
	\]
	where $P_f=\max_{D\in\mathcal D_f}\chi(f,\cO_{\tilde X}(D))$.
	This proves (2).
	
	It remains to prove (3). Let $\cF$ be a rank $r$ vector bundle on $\tilde X$.
	By \eqref{geom:eigensheaves} and the projection formula, we have
	\[
	\chi(f_e,\tilde\pi_e^*\cF)
	=\sum_{j,k=0}^{e-1}\chi(f,\cF(D_{j,k,e})).
	\]
	Let us note that
	\begin{equation*}
		\chi(f,\cF(D_{j,k,e}))-\chi(f,\cF)
		=-\left( c_1(f,\cF)+\frac{r}{2}(D_{j,k,e}-K_{\tilde X}) \right).D_{j,k,e}
	\end{equation*}
	and 
	\begin{equation*}
		\chi(f,\cO_{\tilde X}(D_{j,k,e}))-\chi(f,\cO_{\tilde X})
		=-\frac{1}{2}(D_{j,k,e}-K_{\tilde X}) .D_{j,k,e}.
	\end{equation*}
	So
	\begin{equation}\label{geom:euler-cover}
		\chi(f_e,\tilde\pi_e^*\cF)-r\chi (f_e, \cO_{\tilde X_e})
		=d_e\bigl(\chi(f,\cF)-r\chi (f, \cO_{\tilde X})\bigr)-c_1(f,\cF).S_e,
	\end{equation}
	where $S_e=\sum_{j,k=0}^{e-1}D_{j,k,e}$.
	For each \(i\) the homomorphism
	\((\mathbb Z/e)^2\to \mathbb Z/e,\)
	$(j,k)\to ja_i+kb_i$
	is surjective and every residue occurs \(d_e/e\) times. Hence
	\[
	\sum_{j,k=0}^{e-1}
	\left\lfloor\frac{ja_i+kb_i}e\right\rfloor
	=
	\frac{d_e(e-1)}{2e}(a_i+b_i)
	-\frac{d_e(e-1)}{2e}
	=\frac{d_e(e-1)}{2e}
	\left(a_i+b_i-1\right),
	\]
	and we have
	\[
	S_e=\frac{d_e(e-1)}{2e}
	\sum_i\left(a_i+b_i-1\right)E_i.
	\]
	On the other hand, 
	\[
	\tilde\pi_{e*}R_e
	=\frac{d_e(e-1)}e\left(\sum_i E_i+H_1+H_2\right),
	\]
	so
	\begin{equation*}
		2S_e+\tilde\pi_{e*}R_e
		=\frac{d_e(e-1)}e\bigl(\operatorname{div}(h_1)+\operatorname{div}(h_2)\bigr)
		\sim0.
	\end{equation*}
	So \eqref{geom:euler-cover}, together with the definition of $a$ and \cite[Proposition 3.2]{LaI}, gives
	\[
	a(f_e,\tilde\pi_e^*\cF)
	=d_ea(f,\cF)
	-\frac12c_1(f,\cF)\bigl(2S_e+\tilde\pi_{e*}R_e\bigr)
	=d_ea(f,\cF)
	\]
	as required.
\end{proof}

\begin{Corollary}\label{geom:small-discriminants}
Assume that $k$ has characteristic zero. Let $T$ be the exceptional
$\QQ$-divisor determined by $TE_i=-E_i^2$ for every $i$, and use the
coverings from Proposition~\ref{geom:coverings}. For a reflexive sheaf
$\cE$ of rank $r$ on $X$, let $\cF_e$ be a minimal extension of
$\pi_e^{[*]}\cE$. Then
$0\leq\Delta(f_e,\cF_e)
\leq-\varphi(r)T^2.$ 
In particular,  we have
\[
0\leq
\frac{\Delta(f_e,\cF_e)}{\deg\pi_e}
\leq-\frac{\varphi(r)T^2}{e^2}\longrightarrow0,
\]
uniformly over reflexive sheaves of fixed rank.
\end{Corollary}

\begin{proof}
Put $T_e=\tilde\pi_e^*T/e$. If $C$ lies above $E_i$, then
Proposition~\ref{geom:coverings}(1) and disjointness of the components
above $E_i$ give $T_eC=-C^2>0$. Lemmas~\ref{geom:minimal-slopes}
and~\ref{geom:discriminant-spread} therefore imply
$$\Delta(f_e,\cF_e)\leq-\varphi(r)T_e^2.$$
The lower bound is the relative Bogomolov inequality, and
$T_e^2=(\deg\pi_e)T^2/e^2$ gives the assertion.
\end{proof}

\subsection{Frobenius pullback of a reflexive module}

In this subsection $k$ is an algebraically closed field of characteristic
$p>0$.

The proof below uses the characterization of the relative second Chern
class by Frobenius pullbacks from \cite[Theorem~3.15]{LaI}. A new
ingredient is the following bound for the difference between
a Frobenius pullback of a reflexive module and its reflexive hull:

\begin{Lemma}\label{pos:frobenius-defect}
	Let $R$ be an excellent normal Noetherian local ring of dimension $2$ and characteristic $p>0$.
	There exists a constant $c_R$ such that for every reflexive $R$-module $M$ of rank $r$ and every non-negative
	integer $e$ we have
	$$\operatorname{length}_R\left(
	F_R^{[e]}M/\im(F_R^{e*}M\longrightarrow F_R^{[e]}M)
	\right)\le c_R r p^{2e},$$
	where $F_R^{[e]}M =(F_R^{e*}M)^{**}$.	
\end{Lemma}

\begin{proof}
	Note that $N_e:=F_R^{[e]}M/\im(F_R^{e*}M\to F_R^{[e]}M)$	has finite length as it is supported at the closed point of $\Spec R$.
	
	Let $\mathfrak a_R$ be the ideal consisting of all elements $z\in R$ such that for every reflexive $R$-module $N$ the map $z\,\id_N$ factors through a finite free $R$-module. 
	For $z\in\mathfrak a_R$, choose maps $M\to R^s\to M$ with
	composite $z\,\id_M$. Pulling them back by Frobenius gives
	$$F_R^{e*}M\xrightarrow{\alpha_e}R^s\xrightarrow{\beta_e}F_R^{e*}M,
	\qquad \beta_e\alpha_e=z^{p^e}\,{\id}_{F_R^{e*}M}.$$
	Taking double duals yields
	$$F_R^{[e]}M\xrightarrow{\alpha_e^{**}}R^s
	\xrightarrow{\beta_e^{**}}F_R^{[e]}M$$
	with composite $z^{p^e}\,\id_{F_R^{[e]}M}$. Since $\beta_e^{**}$ is the composite of $\beta_e$ with
	$F_R^{e*}M\to F_R^{[e]}M$, it follows that $z^{p^e}N_e=0$.

	By \cite[Theorem~1.1 and Remark~2.5(1)]{Kimura2025}, $V(\mathfrak a_R)=\operatorname{Sing}R$, so 
	we can choose a system of	parameters $a,b$ contained in $\mathfrak a_R$. 
	Then by the above we have
$$
		\operatorname{length}_R N_e
		\le\operatorname{length}_R
		\bigl(F_R^{[e]}M/(a^{p^e},b^{p^e})F_R^{[e]}M\bigr).
$$
By \cite[Chapter IV.B, Theorem 3]{Serre2000} the length of $F_R^{[e]}M/(a^{p^e},b^{p^e})F_R^{[e]}M$ is equal to 
the Hilbert--Samuel multiplicity $e ({(a^{p^e},b^{p^e})},F_R^{[e]}M)$ of $(a^{p^e},b^{p^e})$ with respect to $F_R^{[e]}M$. But by \cite[Chapter V.A.2]{Serre2000} we have
$$e ({(a^{p^e},b^{p^e})},F_R^{[e]}M)= r e ({(a^{p^e},b^{p^e})},R) =r\operatorname{length}_R\bigl(R/(a^{p^e},b^{p^e})\bigr)
=c_Rrp^{2e},$$
where $c_R=\operatorname{length}_R(R/(a,b))$.
This implies $\operatorname{length}_R N_e	\le c_R rp^{2e}.$
\end{proof}

\subsection{Proof of the linear bound on the correction terms}

\begin{Lemma}\label{geom:exceptional-square}
	Let $D$ and $T$ be exceptional $\mathbb R$-divisors on $\tilde X$.
	If
	\(	|D.E_i|\leq T.E_i	\)
	for every irreducible exceptional curve $E_i$ then
	\(	D^2\geq T^2.	\)
\end{Lemma}

\begin{proof}
	First recall that every $f$-nef exceptional $\mathbb R$-divisor
	$A$ is anti-effective. 	By assumption, both $T+D$ and $T-D$ are $f$-nef.
So $D-T$ is effective and 
$$	D^2-T^2=(D+T).(D-T)\geq0.$$
\end{proof}

\begin{proof}[Proof of Theorem \ref{geom:linear-bound}]
	The upper bound follows from \cite[Theorem 3.12]{LaI}.
	We prove the lower bound.
	
	Assume first that $k$ has characteristic zero.
	For the coverings in Proposition \ref{geom:coverings}, choose
	a minimal extension $\cF_e$ of
	$\pi_e^{[*]}\cE$. Then by Proposition \ref{geom:coverings} (3) we have
	\[
	a(f_e,\cF_e)=\deg\pi_e\cdot a(x,\cE).
	\]
	
	Let $T$ and $T_e=\tilde\pi_e^*T/e$ be as in
	Corollary~\ref{geom:small-discriminants} and its proof. Then 
	$T_eC=-C^2$ for every exceptional component $C$, and
	Corollary~\ref{geom:small-discriminants} gives
	\[
	\Delta(f_e,\cF_e)\leq -\varphi(r)T_e^2.
	\]
Let us set
\[
D_e:=\frac1r c_1(f_e,\cF_e)
-\frac12c_1(f_e,\omega_{\tilde X_e}).
\]
Since
\(	K_{\tilde X_e}C=2g(C)-2-C^2,
	\)
Lemma \ref{geom:minimal-slopes} implies that
	\[
	\left|
D_e.C
	\right|
	\leq-\frac{C^2}{2}=\frac{T_eC}{2}.
	\]
Applying Lemma~\ref{geom:exceptional-square} to $D_e$ and $T_e/2$
	gives
	\[
	D_e^2\geq\frac14T_e^2.
	\]
Hence by definition we have
\begin{align*}
a(f_e,\cF_e)&=\chi(f_e,\cF_e)
	-r\chi (f_e, \cO_{\tilde{X}_e})
	+\frac{r}{2}\left(D_e^2-\frac{1}{4}c_1(f_e,\omega_{\tilde X_e})^2\right)-\frac{\Delta(f_e,\cF_e)}{2r}\\
	&\ge -r\chi (f_e, \cO_{\tilde{X}_e})
	+\left(\frac r8+\frac{\varphi(r)}{2r}\right)T_e^2.
\end{align*}

	By Proposition \ref{geom:coverings} (2), dividing by $\deg\pi_e$
	yields
	\[
	a(x,\cE)\geq-rP_f
	+\left(\frac r8+\frac{\varphi(r)}{2r}\right)\frac{T^2}{e^2}.
	\]
Letting \(e\) tend to infinity gives
	$a(x,\cE)\geq-rP_f$. 
	
	Now assume that $\operatorname{char}k=p>0$, and set
	$\cF=f^{[*]}\cE$. The map $f^*\cE\to\cF$, pulled back by $F_{\tilde X}^e$, induces a map
	$(F_X^{e})^*\cE\to f_*((F_{\tilde X}^{e})^*\cF )$.  So
	we have natural inclusions
	$$\im((F_X^{e})^*\cE\longrightarrow F_X^{[e]}\cE)\subset f_*((F_{\tilde X}^{e})^*\cF)\subset F_X^{[e]}\cE$$
	and by Lemma~\ref{pos:frobenius-defect} there exists some constant $B_x$ such that
	\begin{equation}\label{pos:width-bound}
		\dim H^0\bigl(X,F_X^{[e]}\cE/f_*((F_{\tilde X}^{e})^*\cF)\bigr)
		\le B_x rp^{2e}.
	\end{equation}
	By
	\cite[Lemma~3.8]{LaI}, there is a map
	$\cO_{\tilde X}^{\oplus(r+2)}\to\cF$,
	whose cokernel has finite support. Pulling back  gives
	a map $\cO_{\tilde X}^{\oplus(r+2)}\to(F_{\tilde X}^{e})^*\cF$ with cokernel of
	finite support. This induces a surjection
	$(R^1f_*\cO_{\tilde X})^{\oplus(r+2)}
	\to R^1f_*((F_{\tilde X}^{e})^*\cF).$
	In particular, we have
	\begin{equation}\label{pos:height-bound}
		\dim H^0(X,R^1f_*((F_{\tilde X}^{e})^*\cF))\le(r+2)\dim H^0(X,R^1f_*\cO_{\tilde X}).
	\end{equation}
	Combining \eqref{pos:width-bound} and \eqref{pos:height-bound}, we get
	\begin{equation*}
		0\le\chi(f,(F_{\tilde X}^{e})^*\cF)\le B_x rp^{2e}+(r+2)\dim H^0(X,R^1f_*\cO_{\tilde X}).
	\end{equation*}
	So the Frobenius characterization of the relative second Chern class
	in \cite[Theorem~3.15 and its proof]{LaI} implies  that
	$$-B_x r \le \ch_2(f,\cF)
	=-\lim_{e\to\infty}\frac{\chi(f,F_{\tilde X}^{e*}\cF)}{p^{2e}}\le 0.$$
	Using \cite[Proposition 3.9]{LaI} and the formula
	$$a(x,\cE)=\chi(f,\cF)-r\chi(f,\cO_{\tilde X})+\ch_2(f,\cF)
	-\frac12c_1(f,\cF).K_{\tilde X},$$
	we get the required bound on $|a(x,\cE)|$.
\end{proof}

\section{Properties of relative Chern characters}
\label{sym:section}

We keep the local notation from the preceding section. In particular,
$f:\tilde X\to X$ is a resolution of an algebraic normal surface
singularity, with exceptional divisor $E$.

\subsection{Discriminants of symmetric powers}

In this subsection we prove that the relative discriminants of symmetric powers satisfy
the usual formula holding for discriminants of symmetric powers of a  vector bundle. 
In characteristic zero we use extensions with small
normalized discriminants, whereas in positive characteristic we use
Frobenius and boundedness of the relative discriminants of reflexive
pullbacks of fixed rank. 

Let us recall that for a rank $r$ vector bundle $\cE$, we have  $$r_m:=\rk\Sym^m\cE=\binom{m+r-1}{r-1}$$
and the discriminant satisfies
\begin{equation}\label{sym:ordinary-identities}
	\frac{\Delta(\Sym^m\cE)}{r_m}
	=\binom{m+r}{r+1}\frac{\Delta(\cE)}{r}.
\end{equation}
Indeed, if $x_1,\ldots,x_r$ are the Chern roots of $\cE$, then
\(
\Delta(\cE)=-\sum_{i<j}(x_i-x_j)^2.
\)
The roots of $\Sym^m\cE$ are $\sum_i\alpha_i x_i$, indexed by
the non-negative integers $\alpha_i$ satisfying $\sum_i\alpha_i=m$.
The elementary identity
\[
\sum_{|\alpha|=m}
\left(\sum_i\alpha_i x_i-\frac mr\sum_i x_i\right)^2
=\frac {\binom{m+r}{r+1}}r\sum_{i<j}(x_i-x_j)^2
\]
therefore gives \eqref{sym:ordinary-identities}.

\begin{Lemma}\label{sym:defect}
	For a rank $r$ vector bundle $\cF$ on $\tilde X$, set
	\begin{equation*}
		D_m(f,\cF)
		=\binom{m+r}{r+1}\frac{\Delta(f,\cF)}{r}
		-\frac{\Delta(f,\Sym^m\cF)}{r_m}.
	\end{equation*}
	Then the following assertions hold:
	\begin{enumerate}
		\item We have
		\(
	0\le	D_m(f,\cF)\leq\binom{m+r}{r+1} \frac{\Delta(f,\cF)}{r}.
		\)
		\item If $\tilde\pi:\tilde Y\to\tilde X$ is a proper generically
		finite morphism from a regular integral surface and
		$\tilde Y\xrightarrow{g}Y\xrightarrow{\pi} X$ is the Stein factorization of
		$f\tilde\pi$, then
		\[
		D_m(g,\tilde\pi^*\cF)
		=\deg\tilde\pi\cdot D_m(f,\cF).
		\]
		\item If $\cF$ and $\cF'$ are vector bundles on $\tilde X$ whose
		restrictions to $\tilde X\backslash E$ are isomorphic, then
		\[
		D_m(f,\cF)=D_m(f,\cF').
		\]
	\end{enumerate}
\end{Lemma}

\begin{proof}
The relative Bogomolov	inequality
	\(\Delta(f,\Sym^m\cF)\geq0 \) implies that
$$	D_m(f,\cF)\leq\binom{m+r}{r+1} \frac{\Delta(f,\cF)}{r}.$$	
To prove that $D_m(f,\cF)\geq0$, note that by definition
\[
\Delta(f,\cF)
=\inf\left\{-\frac{\sum_{i<j}(L_i-L_j)^2}{\deg \tilde \pi}\right\},
\]
where the infimum runs over proper generically finite morphisms
$\tilde\pi:\tilde Y\to\tilde X$  and filtrations of
$\tilde\pi^*\cF$ by subbundles with line bundle quotients $\cL_i$.
Here $L_i$ is the relative divisor associated with $\cL_i$.

The induced filtration of $\Sym^m\tilde\pi^*\cF$ has line bundle
quotients $\bigotimes_i\cL_i^{\otimes\alpha_i}$, indexed by
$|\alpha|=m$. Applying \eqref{sym:ordinary-identities} to their
formal roots gives
\[
\frac{\Delta(f,\Sym^m\cF)}{r_m}
\leq-\frac{\binom{m+r}{r+1}}{r}\frac{\sum_{i<j}(L_i-L_j)^2}{\deg \tilde \pi}.
\]
Taking the infimum proves
\[
\frac{\Delta(f,\Sym^m\cF)}{r_m}
\leq \binom{m+r}{r+1}\frac{\Delta(f,\cF)}r,
\]
and hence $D_m(f,\cF)\geq0$.	
	
	Assertion (2) follows from the pullback property of relative
	discriminants and the canonical isomorphism
	\(
	\tilde\pi^*\Sym^m\cF\simeq\Sym^m\tilde\pi^*\cF.
	\)
	
	We prove (3). The bundles, their isomorphism off $E$, and the
	resolution descend from the henselization to an \'etale neighbourhood
	of the algebraic germ. After shrinking this neighbourhood, we
	compactify the normal surface and resolve the boundary without
	changing a neighbourhood of $x$. We thus obtain a proper normal
	surface $\overline X$, smooth away from $x$, and a resolution
	$\overline f:\overline{\tilde X}\to\overline X$ whose exceptional
	locus is $E$.
	
	Extend $\cF$ to a vector bundle $\overline{\cF}$ on
	$\overline{\tilde X}$ by taking a coherent extension followed by its
	double dual. Using the given isomorphism off $E$, glue $\cF'$ near
	$E$ to $\overline{\cF}$ away from $E$. This gives a vector bundle
	$\overline{\cF}'$ agreeing with $\overline{\cF}$ off $E$.
	Let $\overline{\cE}$ be their common reflexive direct image on
	$\overline X$. Then
	\[
	\int_{\overline{\tilde X}}\Delta(\overline{\cF})
	=\int_{\overline X}\Delta(\overline{\cE})
	+\Delta(f,\cF).
	\]
	The analogous identity applies to $\Sym^m\overline{\cF}$, whose
	reflexive direct image is $\Sym^{[m]}\overline{\cE}$.
	Consequently, \eqref{sym:ordinary-identities} gives
	\begin{equation*}
		D_m(f,\cF)
		=\frac{1}{r_m}\int_{\overline X}
		\Delta(\Sym^{[m]}\overline{\cE})
		-\frac{\binom{m+r}{r+1}}{r}\int_{\overline X}\Delta(\overline{\cE}).
	\end{equation*}
	The same calculation for $\overline{\cF}'$ gives the same right
	hand side, proving (3).
\end{proof}

\begin{Theorem}\label{sym:relative-formula}
	For every rank $r$ vector bundle $\cF$ on $\tilde X$ and every
	$m\geq1$, we have
	\[
	\frac{\Delta(f,\Sym^m\cF)}{r_m}
	=\binom{m+r}{r+1}\frac{\Delta(f,\cF)}{r}.
	\]
\end{Theorem}

\begin{proof}
	Put $\cE=f_{[*]}\cF$ and keep $m$ and $r$ fixed.
	First assume that the base field has characteristic zero.
	Using Proposition~\ref{geom:coverings}, choose finite covers
	$\tilde\pi_e:\tilde X_e\to\tilde X$ of degrees $d_e$, with
	induced resolutions $f_e:\tilde X_e\to X_e$, and minimal
	extensions $\cF_e$ isomorphic to $\tilde\pi_e^*\cF$ away
	from the exceptional locus. With $T$ as in
	Corollary~\ref{geom:small-discriminants}, that corollary gives
	\[
	\frac{\Delta(f_e,\cF_e)}{d_e}
	\leq-\frac{\varphi(r)T^2}{e^2}.
	\]
	By Lemma~\ref{sym:defect},
	\[
	0\leq D_m(f,\cF)
	=\frac{D_m(f_e,\cF_e)}{d_e}
	\leq-\binom{m+r}{r+1}\frac{\varphi(r)T^2}{re^2}.
	\]
	Letting $e$ tend to infinity proves $D_m(f,\cF)=0$.
	
	Now assume that the base field has characteristic $p>0$.
	For every $n\geq0$, set
	\(
	\cF_n=f^{[*]}F_X^{[n]}\cE.
	\)
By \cite[Proposition 3.9 and the proof of Theorem 3.12]{LaI},
there is a constant $M\geq0$ such that
	\(\Delta(f,\cF_n)\leq M.	\)
So by Lemma~\ref{sym:defect} we have
	\[
	0\leq D_m(f,\cF)
	=\frac{D_m(f,\cF_n)}{p^{2n}}
	\leq\binom{m+r}{r+1}\frac{\Delta(f,\cF_n)}{rp^{2n}}\le \binom{m+r}{r+1}\frac{M}{rp^{2n}}.
	\]
	Letting $n$ tend to infinity proves $D_m(f,\cF)=0$,
	as required.
\end{proof}

\subsection{Chern characters of direct sums and tensor products}

\begin{Lemma}\label{sym:local-operations}
	Let $f:\tilde X\to X$ be a good resolution of an algebraic normal
	surface singularity with reduced exceptional divisor $E$, and let
	$\cF,\cG$ be vector bundles of ranks
	$r,s$ on $\tilde X$. Then
	\begin{align*}
		\ch_2(f,\cF\oplus\cG)
		&=\ch_2(f,\cF)+\ch_2(f,\cG),\\
		\ch_2(f,\cF\otimes\cG)
		&=s\ch_2(f,\cF)+r\ch_2(f,\cG)
		+c_1(f,\cF)c_1(f,\cG).
	\end{align*}
\end{Lemma}

\begin{proof}
	The definition of the relative second Chern class gives
	$$\ch_2(f,\cF\oplus\cG)
	-\ch_2(f,\cF)-\ch_2(f,\cG)\geq0.$$
	By additivity of the relative Euler characteristic and of the
	first Chern class, this difference equals
	$$a(f,\cF\oplus\cG)-a(f,\cF)-a(f,\cG).$$
	It is therefore unchanged when either bundle is replaced by another
	extension of its restriction to $\tilde X\backslash E$.
	Moreover, it is multiplied by the degree under pullback, by
	\cite[Proposition 3.2]{LaI}.
	
	Assume first that the base field has characteristic zero. Use the
	covers and the divisors $T,T_e$ from
	Corollary~\ref{geom:small-discriminants} and its proof, and choose minimal
	extensions $\cF_e,\cG_e$ of the reflexive pullbacks of
	$f_{[*]}\cF,f_{[*]}\cG$, respectively. For every exceptional
	component $C$, Lemma~\ref{geom:minimal-slopes} gives
	$$\left|\left(\frac{c_1(f_e,\cF_e)}r
	-\frac{c_1(f_e,\cG_e)}s\right)C\right|
	\leq-C^2=T_eC.$$
So applying Lemma~\ref{geom:exceptional-square} to
${c_1(f_e,\cF_e)}/r-{c_1(f_e,\cG_e)}/s$ and $T_e$
gives
\[
\left(
\frac{c_1(f_e,\cF_e)}r-\frac{c_1(f_e,\cG_e)}s
\right)^2\geq T_e^2.
\]
Using the relative Bogomolov inequality for $\cF_e\oplus\cG_e$
	and Corollary~\ref{geom:small-discriminants}, we obtain
	\begin{align*}
		0\leq{}&\ch_2(f_e,\cF_e\oplus\cG_e)
		-\ch_2(f_e,\cF_e)-\ch_2(f_e,\cG_e)\\
		\leq{}&-\frac{rs}{2(r+s)}
		\left(\frac{c_1(f_e,\cF_e)}r
		-\frac{c_1(f_e,\cG_e)}s\right)^2
		+\frac{\Delta(f_e,\cF_e)}{2r}
		+\frac{\Delta(f_e,\cG_e)}{2s}\\
		\leq{}&-\left(\frac{rs}{2(r+s)}
		+\frac{\varphi(r)}{2r}+\frac{\varphi(s)}{2s}\right)T_e^2.
	\end{align*}
	By the invariance and pullback properties above, division by
	$\deg\pi_e$ bounds
	$\ch_2(f,\cF\oplus\cG)-\ch_2(f,\cF)-\ch_2(f,\cG)$
	by a constant times $-T^2/e^2$. Letting $e$ tend to infinity
	proves the first equality.
	In positive characteristic this identity follows from
	\cite[Corollary 3.17]{LaI}.
	
	The case $m=2$ of Theorem~\ref{sym:relative-formula} gives,
	for every rank $t$ vector bundle $\cH$,
	$$\ch_2(f,\Sym^2\cH)
	=(t+2)\ch_2(f,\cH)+\frac12c_1(f,\cH)^2.$$
	Let us apply this to $\cF\oplus\cG$, and use the first equality for the decomposition
	$$\Sym^2(\cF\oplus\cG)
	\simeq\Sym^2\cF\oplus(\cF\otimes\cG)\oplus\Sym^2\cG.$$
	Subtracting the corresponding identities for $\cF$ and $\cG$
	proves the second equality.
\end{proof}

\subsection{Relative asymptotic Riemann--Roch formula}

\begin{Theorem}
	For every vector bundle $\cF$ of rank $r$ on $\tilde X$ one has
	$$\chi(f,\Sym^m\cF)
	=-\frac{m^{r+1}}{(r+1)!}
	\left(c_1(f,\cF)^2-c_2(f,\cF)\right)+O(m^r).$$
\end{Theorem}
\begin{proof}
	Since
	\[
	c_1(f,\Sym^m\cF)=\frac{mr_m}{r}c_1(f,\cF),
	\]
	the definition of $a$ and Theorem~\ref{sym:relative-formula} give
	\begin{equation*}
		\begin{aligned}
			\chi(f,\Sym^m\cF)
			={}&\frac{1}{2r}\binom{m+r}{r+1}\Delta(f,\cF)
			-\frac{m^2r_m}{2r^2}c_1(f,\cF)^2+\frac{mr_m}{2r}c_1(f,\cF)K_{\tilde X}\\
			&
			+r_m\chi(f,\cO_{\tilde X})+a(f,\Sym^m\cF).
		\end{aligned}
	\end{equation*}
	By Theorem~\ref{geom:linear-bound} we have
	\(
	|a(f,\Sym^m\cF)|\leq C_xr_m=O(m^{r-1}).
	\)
	So using
	\[
	\binom{m+r}{r+1}
	=\frac{m^{r+1}}{(r+1)!}+O(m^r),
	\qquad
	r_m=\frac{m^{r-1}}{(r-1)!}+O(m^{r-2}),
	\]
	we obtain
	\begin{align*}
		\chi(f,\Sym^m\cF)
		&=\frac{m^{r+1}}{2r(r+1)!}
		\left(\Delta(f,\cF)-(r+1)c_1(f,\cF)^2\right)
		+O(m^r)\\
		&	=-\frac{m^{r+1}}{(r+1)!}
		\left(c_1(f,\cF)^2-c_2(f,\cF)\right)+O(m^r).
	\end{align*}
\end{proof}

\section{Second Chern classes }

\subsection{The Chern character of reflexive sheaves on surfaces}

Let $X$ be a normal proper surface over an algebraically
closed field. As in \cite[Section 5.5]{LaI}, let
$K^{\refl}(X)$ denote the Grothendieck group of coherent reflexive
sheaves, with relations coming from locally split short exact
sequences. Reflexive tensor product makes $K^{\refl}(X)$
a ring with unit $[\cO_X]$.
We equip $A^*(X)\otimes\RR$ with the ring structure induced
by Mumford's intersection product.

\begin{Theorem}\label{surface:chern-ring}
The Chern character
$\ch:K^{\refl}(X)\to A^*(X)\otimes\RR$ given by
$[\cE]\mapsto\rk\cE+c_1(\cE)+\ch_2(\cE)$
is a homomorphism of rings.
\end{Theorem}

\begin{proof}
Let \(\mathcal E,\mathcal F\) be reflexive sheaves of ranks \(r,s\), respectively, and choose a good resolution \(f:Y\to X\). In $A_0(X)\otimes\RR$,
the definitions give
$$\ch_2(\cE)=f_*\ch_2(f^{[*]}\cE)
-\sum_{x\in\Sing X}\ch_2(f_x,f^{[*]}\cE)[x],$$
where $f_x$ denotes the induced resolution of the local
singularity at $x$. We also have
$$f_*(c_1(f^{[*]}\cE)c_1(f^{[*]}\cF))
=c_1(\cE)c_1(\cF)
+\sum_{x\in\Sing X}
c_1(f_x,f^{[*]}\cE)c_1(f_x,f^{[*]}\cF)[x].$$
Apply the same formulas to $\cF$ and $\cE\hat\otimes\cF$,
using $f^{[*]}\cF$ and $f^{[*]}\cE\otimes f^{[*]}\cF$ as extensions. Ordinary
multiplicativity on $Y$ and Lemma \ref{sym:local-operations} then give
$$\ch_2(\cE\hat\otimes\cF)
=s\ch_2(\cE)+r\ch_2(\cF)+c_1(\cE)c_1(\cF)$$
in $A_0(X)\otimes\RR$.

Reflexive pullback of a locally split short exact sequence is a
locally split short exact sequence of vector bundles on $Y$.
Ordinary additivity on $Y$, together with Lemma \ref{sym:local-operations},
proves additivity of $\ch_2$ for these sequences.
The corresponding identities for rank and first Chern class
are immediate. Thus $\ch$ is well defined, additive and
multiplicative, and $\ch([\cO_X])=1$.
\end{proof}

The following properties of Chern classes of reflexive sheaves on surfaces follow from the results of Section \ref{sym:section}.

\begin{Corollary}\label{sym:surface-formula} Let $\cE$ be a rank $r$ reflexive sheaf  
	on $X$. Then for every $m\geq1$, we have
	$$	\frac{\Delta(\Sym^{[m]}\cE)}{r_m}
	=\binom{m+r}{r+1}\frac{\Delta(\cE)}{r}$$
	in \(A_0(X)\otimes\mathbb R\).
\end{Corollary}

\begin{proof}
	Choose a resolution $f:\tilde X\to X$ and a vector bundle $\cF$
	with $f_{[*]}\cF\simeq\cE$. Let $E$ be the exceptional locus of $f$.
	The definitions of the surface Chern classes  give
	\[
	\Delta(\cE)
	=f_*\Delta(\cF)
	-\sum_{x\in f(E)}\Delta(f_x,\cF)[x]\]
	in $A_0(X)\otimes\RR.$
	Since $f_{[*]}\Sym^m\cF\simeq\Sym^{[m]}\cE$, the same identity,
	together with \eqref{sym:ordinary-identities} and
	Theorem~\ref{sym:relative-formula}, yields
	\begin{align*}
		\frac{\Delta(\Sym^{[m]}\cE)}{r_m}
		&=f_*\left(\frac{\Delta(\Sym^m\cF)}{r_m}\right)
		-\sum_{x\in f(E)}
		\frac{\Delta(f_x,\Sym^m\cF)}{r_m}[x]\\
		&=\frac1r\binom{m+r}{r+1}
		\left(f_*\Delta(\cF)
		-\sum_{x\in f(E)}\Delta(f_x,\cF)[x]\right)=\binom{m+r}{r+1}\frac{\Delta(\cE)}r.
	\end{align*}
\end{proof}

In the following it is convenient to use the \emph{second Segre class} of a reflexive sheaf $\cE$
defined by
$$ s_2(\cE):=c_1(\cE)^2-c_2(\cE)$$
in $A^2(X)\otimes \RR$.

\begin{Corollary}\label{sym:global-rr}
	Let $\cE$ be a rank $r$ reflexive sheaf  
	on $X$. Then we have	
	\begin{align*}
	\chi(X,\Sym^{[m]}\cE)
	=\frac{\int_Xs_2(\cE)}{(r+1)!}m^{r+1}+\frac{r\int_Xs_2(\cE)-c_1(\cE)(c_1(\cE)+K_X)}{2r!}m^r+O(m^{r-1}).
	\end{align*}
\end{Corollary}

\begin{proof}
	By Corollary~\ref{sym:surface-formula} we have
\[
\begin{aligned}
	\operatorname{ch}_2(\operatorname{Sym}^{[m]}\mathcal E)
	={}
	\frac12\left[
	\binom{m+r}{r+1}+\binom{m+r-1}{r+1}
	\right]c_1(\mathcal E)^2-\binom{m+r}{r+1}c_2(\mathcal E).
\end{aligned}
\]
Since $c_1(\Sym^{[m]}\cE)=\binom{m+r-1}{r}c_1(\cE)$, the Riemann--Roch formula and 	Theorem~\ref{geom:linear-bound} imply the required equality.
\end{proof}

\subsection{Second Chern classes in higher dimensions}

Let $X$ be a normal projective variety of dimension $n\geq2$ over
an algebraically closed field $k$. Below $K(X)$ denotes the Grothendieck
group of coherent sheaves, $K_d(X)$ is its subgroup generated by
sheaves with support of dimension at most $d$, and
$$\cF\mapsto c_1(L)\cdot \cF=\cF-L^{-1}\otimes \cF$$
is an additive endomorphism of $K(X)$. 

The following standard lemma is implicitly contained in \cite{LaI}.

\begin{Lemma}\label{high:generic-surface}
Assume that $n>2$. Let $L_1,...,L_{n-2}$ be very ample line bundles on $X$, and let
$S_K$ be the geometric generic complete intersection in
$|L_1|\times...\times|L_{n-2}|$, where $K$ is an algebraic closure
of the function field of this parameter space. 
Then $S_K$ and all the successive intersections are normal and
integral. For every coherent reflexive sheaf $\cE$ on $X$, its
restriction to each of these intersections is reflexive, and
\[
 (\Sym^{[m]}\cE)_{S_K}
 =\Sym^{[m]}(\cE_{S_K})\qquad(m\geq0).
\]
Moreover,
\[
 \chi\left(X,c_1(L_1)...c_1(L_{n-2})\cdot\cE\right)
 =\chi(S_K,\cE_{S_K}).
\]
\end{Lemma}

\begin{Proposition}\label{high:limit}
Let $\cE$ be a reflexive sheaf of rank $r$ on $X$. For arbitrary
line bundles $L_1,...,L_{n-2}$ there exists $a\in\RR$ such that
\[
 \chi\left(X,c_1(L_1)...c_1(L_{n-2})
       \cdot\Sym^{[m]}\cE\right)
 =a m^{r+1}+O(m^r).
\]
\end{Proposition}
\begin{proof}
The computation in \cite[Lemmas 1.6 and 1.7]{LaI} gives
\begin{equation*}
 [\Sym^{[m]}\cE]-\binom{m+r-1}{r-1} [\cO_X]
 -\binom{m+r-1}{r}\bigl([\cE]-r[\cO_X]\bigr)\in K_{n-2}(X).
\end{equation*}
Since every divisor operator lowers the support filtration by one,
applying any $n-1$ such operators, for fixed $M_1,...,M_{n-1}$ we have
\begin{equation}\label{high:lower-dimensional-error}
 \chi\left(X,c_1(M_1)...c_1(M_{n-1})\cdot\Sym^{[m]}\cE\right)=O(m^r).
\end{equation}
Every line bundle $L$ can be written as $A\otimes B^{-1}$ with
$A$ and $B$ very ample. The identity
$$c_1(L)=c_1(A)-c_1(B)+c_1(L)c_1(B)$$
and~\eqref{high:lower-dimensional-error} express the Euler
characteristic in the proposition as the difference of the
corresponding expressions with $L$ replaced by $A$ and $B$,
up to an $O(m^r)$ error. Applying this successively to each $L_i$,
it is enough to prove the assertion when all the $L_i$ are very ample.
In this case Lemma~\ref{high:generic-surface} gives
\[
 \chi\left(X,c_1(L_1)...c_1(L_{n-2})\cdot\Sym^{[m]}\cE\right)
 =\chi\left(S_K,\Sym^{[m]}(\cE_{S_K})\right),
\]
and the assertion follows from Corollary~\ref{sym:global-rr}.
\end{proof}

\begin{Definition}\label{high:definition}
For $\cE\in\Refl(\cO_X)$ and line bundles $L_1,...,L_{n-2}$ we set 
$$\int_Xs_2(\cE)L_1...L_{n-2}
:=(r+1)!\lim_{m\to\infty}
\frac{\chi\left(X,c_1(L_1)...c_1(L_{n-2})
	\cdot\Sym^{[m]}\cE\right)}{m^{r+1}}$$
and
$$\int_Xc_2(\cE)L_1...L_{n-2}
 :=c_1(\cE)^2.L_1...L_{n-2}
   -\int_Xs_2(\cE)L_1...L_{n-2},$$
where \(r=\operatorname{rk}\mathcal E\).
We also set
$$\int_X\ch_2(\cE)L_1...L_{n-2}
 :=\frac12c_1(\cE)^2.L_1...L_{n-2}
   -\int_Xc_2(\cE)L_1...L_{n-2},$$
and define $\int_X\Delta(\cE)=2r \int_Xc_2(\cE)-(r-1)c_1(\cE)^2$.
\end{Definition}

Proposition~\ref{high:limit} and the above definition imply that
	\begin{align*}
		\chi\left(X,c_1(L_1)...c_1(L_{n-2})\cdot\Sym^{[m]}\cE\right)=\frac{m^{r+1}}{(r+1)!}
		\int_Xs_2(\cE)L_1...L_{n-2}+O(m^r).
	\end{align*}

\begin{Theorem}\label{high:properties}
For each reflexive sheaf $\cE$ on $X$, 
$\int_Xc_2(\cE):N^1(X)_{\RR}^{\times(n-2)}\to \RR$ is
a symmetric multilinear form  with the following properties:
\begin{enumerate}
\item If $X$ is a surface, it is the degree of the second Chern class
defined in \cite[Section 4.1]{LaI}.
\item If $\cE$ is locally free, then
$$\int_Xc_2(\cE)L_1...L_{n-2}
 =\int_Xc_2(\cE)\cap c_1(L_1)\cap...\cap c_1(L_{n-2})\cap[X],$$
where the right hand side uses the ordinary Chern class of a vector bundle.
\item If $k\subset K$ is an algebraically closed field extension, then
$$\int_{X_K}c_2(\cE_K)(L_1)_K...(L_{n-2})_K
 =\int_Xc_2(\cE)L_1...L_{n-2}.$$
\item If $n>2$ and $L_1$ is very ample, then for a very general
$H\in|L_1|$,
$$\int_Xc_2(\cE)L_1...L_{n-2}
 =\int_Hc_2(\cE|_H)(L_2)|_H...(L_{n-2})|_H.$$
In particular, if all $L_i$ are very ample, then for a very general
complete intersection surface $S$,
$$\int_Xc_2(\cE)L_1...L_{n-2}=\int_Sc_2(\cE|_S).$$
\end{enumerate}
These forms are uniquely determined by properties (1), (3) and (4).
\end{Theorem}

\begin{proof}
	Proposition~\ref{high:limit} and the identity
	$$c_1(L\otimes M)=c_1(L)+c_1(M)-c_1(L)c_1(M),$$
	together with~\eqref{high:lower-dimensional-error}, show that the
	normalized limit is additive in each line-bundle argument.
	It is symmetric because the divisor operators commute.
	Singular Riemann--Roch shows that the Euler-characteristic
	expressions depend only on the numerical classes of the line
	bundles. Together with the corresponding properties of
	Weil-divisor intersections, this gives the asserted symmetric
	multilinear form on $N^1(X)_{\RR}$.
	
	Property (1) follows from Corollary~\ref{sym:global-rr}.
	Property (3) follows from flat base change for reflexive symmetric
	powers and Euler characteristics, together with invariance of
	Weil-divisor intersections under extension of the ground field.
	Property (4) follows from Lemma~\ref{high:generic-surface},
	after spreading out the reflexive restrictions of the countably
	many symmetric powers and using Definition~\ref{high:definition}.

	We now prove (2). In dimension two, it follows from (1) and
	compatibility of the surface Chern classes with ordinary Chern
	classes of vector bundles.
	In higher dimensions, by multilinearity and (3), we may assume
	that all $L_i$ are very ample and that the ground field is
	uncountable. Choose a very general complete-intersection surface
	$S\in|L_1|\cap\cdots\cap|L_{n-2}|$.
	By (4), the surface case, and the projection formula, we have
	\[
	\begin{aligned}
		\int_Xc_2(\cE)L_1\cdots L_{n-2}
		&=\int_Sc_2(\cE|_S)\cap[S]\\
		&=\int_Xc_2(\cE)\cap c_1(L_1)\cap\cdots
		\cap c_1(L_{n-2})\cap[X],
	\end{aligned}
	\]
	where the right-hand sides use ordinary Chern classes.
	This proves (2).
	
	Finally, properties (1), (3), and (4) determine the values on
	very ample tuples by restriction to a complete-intersection
	surface. Multilinearity then determines all values, since every
	line bundle is the difference of two very ample line bundles.
\end{proof}

\begin{Proposition}\label{high:tensor-product}
	Let $\cE$ and $\cF$ be coherent reflexive sheaves of ranks
	$r$ and $s$ on a normal projective variety $X$ over an
	algebraically closed field. Then
	$$
	\int_X\ch_2(\cE\oplus\cF)=\int_X\ch_2(\cE)+\int_X\ch_2(\cF)
	$$
	and
	\begin{equation*}
		\int_X\ch_2(\cE\hat\otimes\cF)
		=s\int_X\ch_2(\cE)+r\int_X\ch_2(\cF)+c_1(\cE)c_1(\cF).
	\end{equation*}
\end{Proposition}

\begin{proof}
By properties (1), (3) and (4) of Theorem~\ref{high:properties}
and multilinearity, it is enough to prove the identities on a normal
projective surface. In this case the equalities follow from
Theorem~\ref{surface:chern-ring}.
\end{proof}

\begin{Remark}
The above proposition and Corollary \ref{sym:surface-formula} extend to computations of the second Chern character of the reflexive sheaf associated to a given rank $r$ reflexive sheaf by some rational \(\mathrm{GL}_r\)-representation.
\end{Remark}

\subsection{The Riemann--Roch inequality on normal varieties}

\begin{Theorem}\label{high:linear-rr}
Let $X$ be a normal proper algebraic variety of dimension $n\ge 2$ over an algebraically
closed field. Let $L_1,...,L_{n-2}$ be very ample line bundles on $X$. There is
a constant $C_X$ depending only on \(X,L_1,\ldots,L_{n-2}\) such that for every reflexive sheaf $\cE$ of rank $r$
one has
\begin{align*}
\Big|\chi\left(X,c_1(L_1)...c_1(L_{n-2})\cdot\cE\right)
 &-r\chi\left(X,c_1(L_1)...c_1(L_{n-2})\cdot\cO_X\right)
 -\int_X\ch_2(\cE)L_1...L_{n-2}\\
 &+\frac12c_1(\cE).(K_X+L_1+...+L_{n-2}).L_1...L_{n-2}\Big|
 \leq C_Xr.
\end{align*}
\end{Theorem}

\begin{proof}
If $n=2$ then \cite[Theorem 4.4]{LaI} implies that 
	\begin{align*}
	\left|\chi(X,\cE)-r\chi(X,\cO_X) -\int_X\ch_2(\cE) +\frac12c_1(\cE).K_X \right|= \left|\sum_{x\in\Sing X}a(x,\cE)\right|.
\end{align*}
So the assertion follows from 	Theorem~\ref{geom:linear-bound}.

If $n>2$, let $S_K$ be the fixed geometric generic surface from
Lemma~\ref{high:generic-surface}. For every reflexive sheaf
$\cE$ on $X$, that lemma gives reflexivity of $\cE_{S_K}$ and
identifies the Euler characteristics in the assertion with those
of $\cE_{S_K}$ and $\cO_{S_K}$. By
Theorem~\ref{high:properties},
$$\int_X\ch_2(\cE)L_1...L_{n-2}
 =\int_{S_K}\ch_2(\cE_{S_K}).$$
Moreover, the canonical divisor of $S_K$ is the restriction of
$K_X+L_1+...+L_{n-2}$. Thus the expression inside the absolute
value is
$$\sum_{x\in\Sing S_K}a(x,\cE_{S_K}),$$
and the assertion follows from  Theorem~\ref{geom:linear-bound}.
\end{proof}

\subsection{Chern classes of pullbacks by finite morphisms}
\begin{Theorem}\label{thm:c_2-mult-by-degree}
	If $\pi:Y\to X$ is a finite surjective morphism from a normal
	projective variety, then for each reflexive sheaf $\cE$ on $X$ we have
	$$\int_Yc_2(\pi^{[*]}\cE)\pi^*L_1...\pi^*L_{n-2}
	=\deg\pi\cdot\int_Xc_2(\cE)L_1...L_{n-2}$$
for arbitrary line bundles \(L_1,\ldots,L_{n-2}\) on \(X\).
\end{Theorem}

\begin{proof}
	Assume first that $k$ has characteristic zero and all $L_i$ are
	very ample. After extending the base field, let $S$ be the geometric
	generic complete intersection in $|L_1|\times...\times|L_{n-2}|$,
	and put $T=\pi^{-1}(S)$. By Bertini's theorem, $T$ is a normal
	surface. The incidence and Koszul arguments in the proof of
	Lemma~\ref{high:generic-surface} also apply to the pulled-back
	basepoint-free linear systems on $Y$. Thus all the reflexive
	symmetric powers under consideration restrict reflexively, and
	Definition~\ref{high:definition} identifies the required Chern
	numbers with those on $S$ and $T$.
	Moreover,
	$(\pi^{[*]}\cE)|_T\simeq(\pi|_T)^{[*]}(\cE|_S)$,
	because this equality holds away from finitely many points and both
	sides are reflexive. The finite morphism $\pi|_T:T\to S$ has degree
	$\deg\pi$, so \cite[Proposition 4.2]{LaI} implies that
	$$\int_Tc_2((\pi^{[*]}\cE)|_T)
	=\deg\pi\cdot\int_Sc_2(\cE|_S),$$
	as required. Multilinearity gives the assertion for arbitrary line bundles.
	
	Assume now that $k$ has characteristic $p>0$.
	By restriction to a geometric generic complete-intersection surface,
	Theorem~\ref{high:properties} and the surface Frobenius formula
	\cite[Proposition 4.2]{LaI} give
	\begin{equation}\label{high:frobenius-character}
		\int_X\ch_2(F_X^{[e]}\cE)
		=p^{2e}\int_X\ch_2(\cE).
	\end{equation}
	The same equality holds for reflexive sheaves on $Y$.
	
	By multilinearity, it is enough to prove the required equality when
	both $L_i$ and $\pi^*L_i$ are very ample: every line bundle is a
	difference of ample line bundles, and sufficiently large powers of
	an ample line bundle and its finite pullback are very ample.
	Put $d=\deg\pi$, $r=\rk\cE$, and $\cA=\pi_*\cO_Y$.
	For $e\geq0$, put $q=p^e$ and
	$\cB_e=F_X^{[e]}\cE\hat\otimes\cA.$
	We have a natural isomorphism
	$\cB_e\simeq\pi_*F_Y^{[e]}(\pi^{[*]}\cE),$
as  both sides are reflexive and agree on the big open subset
	of $X$ where $\cE$ is locally free.
	
	The sheaf $\cB_e$ has rank $rd$ and
	$c_1(\cB_e)=dq\,c_1(\cE)+r\,c_1(\cA).$
	Moreover, Proposition~\ref{high:tensor-product} and
	\eqref{high:frobenius-character} imply that
	$$\int_X\ch_2(\cB_e)
	=dq^2\int_X\ch_2(\cE)
	+q\,c_1(\cE)c_1(\cA)+r\int_X\ch_2(\cA).$$
	Applying Theorem~\ref{high:linear-rr} to $\cB_e$, we obtain
	\begin{align*}
		\chi\left(X,c_1(L_1)...c_1(L_{n-2})\cdot\cB_e\right)
		&=\int_X\ch_2(\cB_e)L_1...L_{n-2}+O(q)\\
		&=dq^2\int_X\ch_2(\cE)L_1...L_{n-2}+O(q).
	\end{align*}
	Here the rank is fixed and the first-Chern-class term is $O(q)$.
	Similarly, applying Theorem~\ref{high:linear-rr} on $Y$ gives
	\begin{align*}
		\chi\left(Y,c_1(\pi^*L_1)...c_1(\pi^*L_{n-2})
		\cdot F_Y^{[e]}(\pi^{[*]}\cE)\right)=q^2\int_Y\ch_2(\pi^{[*]}\cE)
		\pi^*L_1...\pi^*L_{n-2}+O(q).
	\end{align*}
	Since $\pi$ is finite, the two Euler characteristics agree by the
	projection formula. Dividing by $q^2$ and letting $e\to\infty$, we obtain
	\begin{equation}\label{high:finite-character}
		\int_Y\ch_2(\pi^{[*]}\cE)\pi^*L_1...\pi^*L_{n-2}
		=d\int_X\ch_2(\cE)L_1...L_{n-2}.
	\end{equation}
	
	Finally, the second Chern class of a rank-one reflexive sheaf
	vanishes, by (1), (4) of Theorem~\ref{high:properties},
	multilinearity, and \cite[Proposition 4.2(1)]{LaI}.
	Applying~\eqref{high:finite-character} to $\det\cE$ and using
	$\pi^{[*]}(\det\cE)\simeq\det(\pi^{[*]}\cE)$, we obtain
	$$c_1(\pi^{[*]}\cE)^2.\pi^*L_1...\pi^*L_{n-2}
	=d\,c_1(\cE)^2.L_1...L_{n-2}.$$
	Together with~\eqref{high:finite-character}, this proves the assertion.
\end{proof}

\section{Approximation of discriminants by alterations}

\begin{Definition}
	Let $X$ be a normal projective variety  over an algebraically closed
	field $k$ and let $\cE$ be a reflexive sheaf on $X$. 
	Let  $f:Z\to X$ be an alteration from a smooth
	projective variety $Z$.	We say that a reflexive sheaf $\cF$ on $Z$
	is \emph{compatible with $\cE$} if	\(	g_{[*]}\cF\simeq h^{[*]}\cE,\) 	
	where \( Z\xrightarrow{g}Y\xrightarrow{h}X \) is
	the Stein factorization of $f$.	 In this case we say that $(f, \cF)$ is \emph{compatible with $(X,\cE)$}.
\end{Definition}

\subsection{Approximation of discriminants on surfaces}

\begin{Proposition}\label{app:global-discriminant-approximation}
	Let $X$ be a normal projective surface over an algebraically
	closed field $k$, and fix a projective good resolution
	$f:S\to X$ which is an isomorphism over $X_{\rm reg}$.
	There exists a sequence of alterations
	$\widetilde\pi_i:S_i\to S$ from smooth projective surfaces,
	with degrees $d_i\to\infty$, such that for every positive
	integer $r$ there is a constant $C_{f,r}\geq0$ with the following
	property. Every rank $r$ reflexive sheaf $\cE$ on $X$ admits
	compatible vector bundles $\cF_i$ on $S_i$ satisfying 
	\[
	\cF_i|_{(f\widetilde\pi_i)^{-1}(X_{\rm reg})}
	\simeq
	(f\widetilde\pi_i)^*(\cE|_{X_{\rm reg}})
	\]
	and
	\[
	0\leq
	\frac{1}{d_i}\int_{S_i}\Delta(\cF_i)
	-\int_X\Delta(\cE)
	\leq\frac{C_{f,r}}{i^2}.
	\]
	In particular,
	\[
	\int_X\Delta(\cE)
	=\lim_{i\to\infty}
	\frac{1}{d_i}\int_{S_i}\Delta(\cF_i).
	\]
\end{Proposition}

\begin{proof}
	Assume first that $\operatorname{char}k=0$.
	Choose a nonempty finite set $B\subset X$ containing
	$\operatorname{Sing}X$. For $x\in B$, write $X_x=\Spec\cO_{X,x}^h$ and $S_x=S\times_X X_x.$
	Choose an increasing sequence of sufficiently large primes
	$e_i\geq i$. Proposition~\ref{geom:coverings} gives finite
	coverings $S_{x,i}\to S_x$ of degree $e_i^2$, whose composite with $f_x$ has Stein factorization
	\[
	S_{x,i}\xrightarrow{f_{x,i}}X_{x,i}
	\xrightarrow{\pi_{x,i}}X_x.
	\]
	By Corollary~\ref{geom:small-discriminants}, minimal extensions
	$\cG_{x,i}$ of $\pi_{x,i}^{[*]}(\cE|_{X_x})$ satisfy
	\[
	0\leq\Delta(f_{x,i},\cG_{x,i})\leq C_{x,r},
	\]
	where $C_{x,r}$ is independent of $i$ and $\cE$.
	
	We now construct one alteration dominating these local coverings.
	Put $K=K(X)$ and $K_x=K(X_x)$.
	Since $K_x/K$ is algebraic, generators of
	$K(S_{x,i})/K_x$ are algebraic over $K$.
	Choose a finite Galois extension $L_i/K$ containing such
	generators for every $x\in B$.
	Every field factor of $L_i\otimes_K K_x$ then contains a
	$K_x$-isomorphic copy of $K(S_{x,i})$: this holds for one
	factor by construction, and the Galois group acts transitively
	on the factors.
	
	Let $Y_i$ be the normalization of $X$ in $L_i$, and let $S_i$
	be a smooth projective resolution of the normalization of $S$
	in $L_i$. We obtain
	\[
	\widetilde\pi_i:S_i\to S,\qquad
	S_i\xrightarrow{g_i}Y_i\xrightarrow{\pi_i}X.
	\]
	Thus, over the henselization at every $y\in\pi_i^{-1}(x)$,
	the surface $S_i$ dominates $S_{x,i}$.
	Write $b_{i,y}$ for the degree of this local alteration.
	With $d_i=[L_i:K]$, we have
	\[
	\sum_{y\in\pi_i^{-1}(x)}b_{i,y}=\frac{d_i}{e_i^2}.
	\]
	In particular, $d_i\geq e_i^2\to\infty$.
	
	Pull back $\cG_{x,i}$ along these local alterations.
	The resulting bundles glue with
	$(f\widetilde\pi_i)^*\cE$ over
	$(f\widetilde\pi_i)^{-1}(X\setminus B)$ to a vector bundle
	$\cF_i$ on $S_i$.
	Indeed, relative to
	$\widetilde\pi_i^*f^{[*]}\cE$, these are modifications supported
	on the disjoint fibers over $\pi_i^{-1}(B)$.
	Each modification is determined on a finite infinitesimal
	neighborhood of its fiber, and hence descends from the
	henselization. The resulting bundle is compatible with $\cE$.
	
	Relative discriminants multiply by the degree under pullback.
	The other exceptional fibers contribute zero, since there
	$\cF_i$ is the pullback of a vector bundle on $Y_i$.
	Consequently,
	\[
	\begin{aligned}
		0
		\leq \int_{S_i}\Delta(\cF_i)
		-d_i\int_X\Delta(\cE)
		=\sum_{x\in B}\sum_{y\in\pi_i^{-1}(x)}
		b_{i,y}\Delta(f_{x,i},\cG_{x,i})
		\leq\frac{d_i}{e_i^2}\sum_{x\in B}C_{x,r}.
	\end{aligned}
	\]
	Dividing by $d_i$ proves the assertion with
	$C_{f,r}=\sum_{x\in B}C_{x,r}$.
	
	Now assume that $\operatorname{char}k=p>0$.
	Take $S_i=S^{(-i)}$, $Y_i=X^{(-i)}$, and $g_i=f^{(-i)}$,
	where $(-i)$ denotes the inverse Frobenius twist.
	Let $\widetilde\pi_i:S_i\to S$ and $\pi_i:Y_i\to X$
	be the relative Frobenius morphisms, and set
	$\cF_i=g_i^{[*]}\pi_i^{[*]}\cE.$
	By \cite[Proposition 3.9 and the proof of Theorem 3.12]{LaI},
	the sum of the local relative discriminants of $\cF_i$
	is bounded by a constant $C_{f,r}$ independent of $i$
	and $\cE$; these bounds are unchanged by Frobenius twists.
	Hence
	\[
	0\leq
	\frac{1}{p^{2i}}\int_{S_i}\Delta(\cF_i)
	-\int_X\Delta(\cE)
	\leq\frac{C_{f,r}}{p^{2i}}
	\leq\frac{C_{f,r}}{i^2}.
	\]
\end{proof}

\begin{Remark}
The displayed isomorphism over $X_{\rm reg}$ is stronger than
compatibility, that is, $g_{i[*]}\cF_i\simeq\pi_i^{[*]}\cE$.
It does not imply that
$\cF_i\simeq(f\widetilde\pi_i)^{[*]}\cE$.
\end{Remark}

\begin{Corollary}\label{app:discriminant-infimum}
	Let $X$ be a normal projective surface over an algebraically
	closed field. For every reflexive sheaf $\cE$ on $X$ of positive
	rank, we have
	\[
	\int_X\Delta(\cE)
	=
	\inf_{f,\cF}
	\frac{\int_Z\Delta(\cF)}{\deg f},
	\]
	where $f:Z\to X$ runs through all alterations from smooth projective surfaces, and $\cF$ runs through all
	vector bundles on $Z$ that are  compatible with $\cE$.	
\end{Corollary}

\begin{proof}
	Let $(f,\cF)$ be a pair compatible with $(X, \cE)$ and let \( Z\xrightarrow{g}Y\xrightarrow{h}X \) be
	the Stein factorization of $f$.
	The birational discriminant inequality and finite-pullback
	functoriality give
	\[
	\int_Z\Delta(\cF)\geq	\int_Y\Delta(g_{[*]}\cF)
	=	\int_Y\Delta(h^{[*]}\cE)=(\deg f)	\int_X\Delta(\cE).
	\]
	The reverse inequality follows from the approximating sequence
	in Proposition~\ref{app:global-discriminant-approximation}.
\end{proof}

\begin{Example}\label{need-for-alterations}
	The finite coverings cannot in general be omitted from the infimum.
	Let $X=\PP(1,1,2)$, the projective quadric cone, let
	$\cL=\cO_X(1)$, and put $\cE=\cO_X\oplus\cL$.
	Then $c_1(\cL)^2=1/2$, and hence
	\(
	\int_X	\Delta(\cE)=-1/2.
	\)
	For every vector bundle extension $\cF$ on a smooth resolution of $X$,
	the ordinary discriminant is an integer and it is at least $-1/2$ by the
	local discriminant formula. It is therefore nonnegative.
	On the minimal resolution $\mathbb F_2\to X$, the bundle
	$\cO_{\mathbb F_2}\oplus\cO_{\mathbb F_2}(F)$, where $F$ is a fiber
	of the ruling, extends $\cE$ and has discriminant zero.
	Thus the infimum over resolutions of $X$ itself is zero, whereas the
	normalized infimum allowing finite covers is $-1/2$.
\end{Example}

\subsection{Approximation of discriminants in higher dimensions}

The following lemma generalizes \cite[Proposition 4.2, (3)]{LaI} to higher dimensions.

\begin{Lemma}\label{birational-maps}
	Let $g:Z\to Y$ be a projective birational morphism from a smooth
projective variety to a normal projective variety of dimension $n\ge 2$. Let
$\cF$ be a vector bundle on $Z$. Then for any nef line bundles $H_1,\ldots,H_{n-2}$ on $Y$ we have
\begin{equation*}
	\int_Z\Delta(\cF)g^*H_1\cdots g^*H_{n-2}
	\geq
	\int_Y\Delta(g_{[*]}\cF)H_1\cdots H_{n-2}.
\end{equation*}
\end{Lemma}

\begin{proof}
Let us set $\cG:=g_{[*]}\cF$ and assume that all $H_i$ are ample.
Replacing them by sufficiently large multiples and extending
the ground field if necessary, choose a sufficiently general
	normal complete-intersection surface
	$S\in |H_1|\cap\cdots\cap |H_{n-2}|$
	such that $\cG|_S$ is reflexive and the restriction formula
	for its discriminant holds.
	
	We can also arrange that $S$ avoids the locus over which $g$
	has fibers of dimension at least two.
	Indeed, this locus has codimension at least three: its inverse
	image is contained in the exceptional locus of $g$, which has
	dimension at most $n-1$.
	The nonisomorphism locus of $g$ meets $S$ in finitely many points.
	
	Set $W=g^{-1}(S)$.
	Every irreducible component of $W$ has dimension at least two,
	whereas every fiber over the exceptional points of $S$ has
	dimension at most one.
	Thus $W$ has a unique irreducible component, dominating $S$,
	and has dimension two.
	It is a complete intersection in the smooth variety $Z$, hence
	Cohen--Macaulay, and it is generically reduced because $g$ is
	an isomorphism over a dense open subset of $S$.
	Consequently, $W$ is integral.
	
	Let $\rho:T\to W$ be a resolution of its normalization and put
	$\cV=\rho^*(\cF|_W)$.
	The induced morphism $T\to S$ is birational, and the reflexive
	direct image of $\cV$ is $\cG|_S$.
	The projection formula and the surface discriminant inequality
	therefore give
	\[
	\begin{aligned}
		\int_Z\Delta(\cF)g^*H_1\cdots g^*H_{n-2}
		&=\int_W\Delta(\cF|_W)=\int_T\Delta(\cV)\\
		&\geq\int_S\Delta(\cG|_S)=\int_Y\Delta(\cG)H_1\cdots H_{n-2}.
	\end{aligned}
	\]
	This proves the assertion for ample line bundles. The nef case follows by continuity.
\end{proof}

The following theorem is a higher dimensional version of Corollary \ref{app:discriminant-infimum}. 
We use alterations of general complete intersection surfaces coming from this corollary to reduce the higher dimensional version to the surface case.

\begin{Theorem}\label{app:discriminant-infimum-higher}
	Let $X$ be a normal projective variety of dimension $n\geq2$
	over an algebraically closed field, and let
	$L_1,\ldots,L_{n-2}$ be nef line bundles on $X$.
	For every reflexive sheaf $\cE$ of positive rank on $X$, we have
	\[
	\int_X\Delta(\cE)L_1\cdots L_{n-2}
	=
	\inf_{f,\cF}
	\frac{
		\int_Z\Delta(\cF)f^*L_1\cdots f^*L_{n-2}
	}{\deg f},
	\]
	where $f:Z\to X$ runs through all alterations from smooth
	projective varieties, and $\cF$ runs through all vector bundles
	on $Z$ compatible with $\cE$.
\end{Theorem}

\begin{proof}
	Let $(f,\cF)$ be compatible with $(X,\cE)$, and let
	$Z\xrightarrow{g}Y\xrightarrow{h}X$
	be the Stein factorization of $f$.
	By Lemma~\ref{birational-maps} and Theorem \ref{thm:c_2-mult-by-degree} we have
	\[
	\begin{aligned}
		\int_Z\Delta(\cF)f^*L_1\cdots f^*L_{n-2}
		&\geq
		\int_Y\Delta(h^{[*]}\cE)h^*L_1\cdots h^*L_{n-2}\\
		&=(\deg f)\int_X\Delta(\cE)L_1\cdots L_{n-2}.
	\end{aligned}
	\]
	This proves the required lower bound.
	
	If $n=2$, the opposite inequality follows from
	Corollary~\ref{app:discriminant-infimum}.
	Assume that $n\geq3$.
	By extension of the algebraically closed ground field and
	specialization of the resulting alterations and bundles,
	we may assume that the ground field is uncountable.
	
	Choose a sufficiently ample line bundle $H$ so that
	$H\otimes L_i^{-1}$ is nef for every $i$.
	For each compatible pair, the difference
	\[
	\frac{\int_Z\Delta(\cF)f^*L_1\cdots f^*L_{n-2}}{\deg f}
	-\int_X\Delta(\cE)L_1\cdots L_{n-2}
	\]
	is multilinear in the $L_i$ and nonnegative on nef tuples.
	It therefore does not decrease when any $L_i$ is replaced
	by $H$. Hence it suffices to prove the assertion when
	$L_1=\cdots=L_{n-2}=H$.
	
	Let $U\subset X$ be the open subset where $X$ is smooth
	and $\cE$ is locally free. Its complement has codimension
	at least two.
	Choose $n-1$ very general sections of $H$ whose common zero
	locus is a smooth curve $C\subset U$.
	They define a rational map $X\dashrightarrow\PP^{n-2}$,
	which is resolved by the blowup
	$q:P=\operatorname{Bl}_C X\longrightarrow X.$
	Let $\lambda:P\to\PP^{n-2}$ be the induced morphism and
	$E$ the exceptional divisor. Then
	\begin{equation}\label{app:pencil-divisor}
		q^*H=\lambda^*\cO_{\PP^{n-2}}(1)+E.
	\end{equation}
	
	Let $S$ be the geometric generic fibre of $\lambda$.
	It identifies with the corresponding complete-intersection
	surface in $X$, since $C$ is Cartier on that surface.
	Our choice of sections ensures that $S$ is normal,
	$\cE|_S$ is reflexive, and
	\begin{equation*}
		\int_S\Delta(\cE|_S)
		=
		\int_X\Delta(\cE)H^{n-2}.
	\end{equation*}
	Fix $\varepsilon>0$.
	By Proposition~\ref{app:global-discriminant-approximation},
	there exist an alteration $a_\eta:V\to S$ and a vector
	bundle $\cV$ on $V$ such that
	\begin{equation}\label{app:generic-approximation}
		0\leq
		\frac{\int_V\Delta(\cV)}{\deg a_\eta}
		-\int_X\Delta(\cE)H^{n-2}
		<\varepsilon.
	\end{equation}
	Moreover, $\cV$ agrees with the pullback of $\cE$
	over the inverse image of $U$.
	
	Descending these data to a finite extension of
	$k(\PP^{n-2})$ and taking a projective closure gives a
	projective alteration $w:W\to P$.
	Over a dense open subset of $\PP^{n-2}$, the descended
	bundle and its isomorphism with the pullback of $\cE$
	over $U$ extend.
	Glue this extension with $(qw)^*(\cE|_U)$ on
	$(qw)^{-1}(U)$, and extend the resulting vector bundle
	to a coherent sheaf $\cM$ on $W$.
	
	Let $\beta:W'\to W$ be the blowup of the Fitting ideal
	$\operatorname{Fitt}_{\operatorname{rk}\cE}(\cM)$.
	By \cite[Tag 0ESN]{StacksProject}, the sheaf
	$\cV'=\beta^*\cM/\mathrm{torsion}$
	is locally free, and $\beta$ is an isomorphism wherever
	$\cM$ is locally free.
		By de Jong's alteration theorem \cite{deJong1996}, there
	exists an alteration $t:Z\to W'$ with $Z$ smooth and
	projective. Set
	\[
	a=w\beta t,\qquad f=qa,\qquad \cF=t^*\cV'.
	\]
	Then $(f,\cF)$ is compatible with $(X,\cE)$ and the restriction of $\cF$ to $a^*E$ is pulled back from
	$\cE|_C$. Since $C$ is a smooth curve, this gives $\Delta(\cF)\cdot a^*E=0.$
So using \eqref{app:pencil-divisor} we obtain
	\[
	\begin{aligned}
		\frac{\int_Z\Delta(\cF)(f^*H)^{n-2}}{\deg f}
		=
		\frac{
			\int_Z\Delta(\cF)
			(\lambda a)^*c_1(\cO_{\PP^{n-2}}(1))^{n-2}
		}{\deg f}=\frac{\int_V\Delta(\cV)}{\deg a_\eta}.
	\end{aligned}
	\]
	Here the finite extension used in the descent and the
	alteration $t$ multiply the relevant intersection numbers
	and degrees by the same factors.
	The required upper bound now follows from
	\eqref{app:generic-approximation}.
\end{proof}

\begin{Remark}
	The nef hypothesis in Theorem \ref{app:discriminant-infimum-higher} cannot be replaced by arbitrary Cartier divisors. For example, let $g:Z\to X$ be the blowup of a smooth projective
	threefold $X$ along a smooth curve $C\subset X$, with exceptional
	divisor $E$, and let $L$ be a Cartier divisor on $X$. Set
	$	\mathcal E=\mathcal O_X^{\oplus2}$ and $
	\mathcal F_t=\mathcal O_Z(tE)\oplus\mathcal O_Z(-tE).$
	Then \(g_{[*]}\mathcal F_t\simeq\mathcal E\) but
	\[
	\int_Z\Delta(\mathcal F_t)g^*L
	=4t^2(L\cdot C).
	\]
	If \(L\cdot C<0\), the infimum is \(-\infty\), whereas the discriminant of \(\mathcal E\) is zero.
\end{Remark}

\begin{Corollary}\label{Bogomolov-inequality}
Let $X$ be a normal projective variety of dimension $n\geq2$
over an algebraically closed field, and let
$L=(L_1,\ldots,L_{n-1})$ be a collection of nef  line bundles on $X$ such that the $1$-cycle $L_1\ldots L_{n-1}$ is numerically non-trivial. If $\cE$ is a strongly slope $L$-semistable reflexive sheaf then  
$$ \int_X \Delta (\cE)L_2\dots L_{n-1}\ge 0.$$ 
\end{Corollary}

\begin{proof}
	Let $(f,\cF)$ be a pair compatible with $(X, \cE)$.
 Then $\cF$	is strongly $(f^*L_1,\ldots,f^*L_{n-1})$-semistable and 
 Bogomolov's inequality on $Z$ (see \cite[Theorem 3.4]{Langer2004})
 gives
 \[
 \int_Z\Delta(\cF)f^*L_2\cdots f^*L_{n-1}\geq0.
 \]
 So the inequality follows from Theorem  \ref{app:discriminant-infimum-higher}.
\end{proof}

\subsection{Relation to Chen's results} \label{Chen}
	
Let $X$ be a normal complex projective variety of dimension
	$n\geq2$, and let $L_1,\ldots,L_{n-2}$ be ample line bundles.	
	In this setting,  Chen in \cite{Chen2025} defined the discriminant of a reflexive sheaf $\cE$
as
	\[
	\int_X\Delta_{\mathrm{Chen}}(\cE)L_1\cdots L_{n-2}
	=
	\inf_{f,\cF}
	\int_Z\Delta(\cF)f^*L_1\cdots f^*L_{n-2},
	\]
	where $f:Z\to X$ runs through projective birational morphisms
	from smooth projective varieties, and $\cF$ runs through
	reflexive sheaves satisfying $f_{[*]}\cF\simeq\cE$. His definition used analytic sheaves
	and K\"ahler resolutions but a compact K\"ahler resolution of a projective variety
	is projective, and GAGA identifies its coherent analytic sheaves with algebraic ones. So the above formula gives the invariant defined in \cite{Chen2025}.

Note that the infimum is unchanged if we restrict to vector bundles.
So by Theorem~\ref{app:discriminant-infimum-higher}, we have
	\[
	\int_X\Delta_{\mathrm{Chen}}(\cE)L_1\cdots L_{n-2}\ge 	\int_X\Delta(\cE)L_1\cdots L_{n-2}.
	\]
	This inequality can be strict, as shown by
	Example~\ref{need-for-alterations}.

	Suppose now that $\cE$ is slope semistable with respect to a collection of ample line bundles
	$(L_1,\ldots,L_{n-1})$. In this case Chen also proved a Bogomolov type inequality (see \cite[Theorem 1.14]{Chen2025}). His inequality follows immediately from Bogomolov's inequality on the resolutions of $X$
	(as in the proof of Corollary \ref{Bogomolov-inequality}) and it is usually weaker than that in Corollary \ref{Bogomolov-inequality}. The  non-trivial part of \cite[Theorem 1.14]{Chen2025} is the second one that concerns equality $\int_X\Delta_{\mathrm{Chen}}(\cE)L_2\cdots L_{n-1}=0$.
Then \cite[Theorem 1.14]{Chen2025} gives a filtration with torsion-free quotients whose restrictions to $X_{\rm reg}$
are projectively unitary flat.
Corollary \ref{Bogomolov-inequality} implies that in this case \(\int_X\Delta(\cE)L_2\cdots L_{n-1}=0. \)

We expect the following stronger statement.
If $\cE$ is slope polystable with respect to
$(L_1,\ldots,L_{n-1})$ and
\[
\int_X\Delta(\cE)L_2\cdots L_{n-1}=0,
\]
then $\cE|_{X_{\rm reg}}$ is locally free and projectively unitary flat.
More precisely, writing $r=\rk\cE$, we expect that the analytic projectivization of $\cE|_{X_{\rm reg}}$ is induced by
a representation \(\pi_1(X_{\rm reg}^{\rm an})\to\operatorname{PU}(r).\)
This conclusion does not follow from the comparison above,
since vanishing of our discriminant does not imply
vanishing of Chen's invariant.

Conversely, projective unitary flatness of $\cE|_{X_{\rm reg}}$ implies
the asserted discriminant vanishing. Let us set $\cA=\cE\hat\otimes\cE^*$.
Then $\cA|_{X_{\rm reg}}$ is unitary flat and
\(
\int_X\ch_2(\cA)=-\int_X\Delta(\cE).
\)
So applying \cite[Lemma 6.12]{Langer2026WP} to the restriction of $\cA$
to a very general normal complete-intersection surface gives
\[\int_X\Delta(\cE)L_2\cdots L_{n-1}=-
\int_X\ch_2(\cA)L_2\cdots L_{n-1}=0.
\]
For more details and an interpretation of vanishing of the discriminant in the semistable case, we refer to \cite[Sections 5 and 6]{Langer2026WP}.

\section{Applications}\label{Applications}

In this section we prove strong restriction theorems for semistable sheaves on normal varieties
and we show some new boundedness results in characteristic zero.  We also prove a Miyaoka type inequality and some hyperbolicity results on normal surfaces.

\medskip

Let $X$ be a normal projective variety of dimension $n$ defined over an algebraically closed field $k$ and let $L=(L_1,...,L_{n-1})$ be a collection of nef line bundles on $X$ such that the $1$-cycle $L_1...L_{n-1}$ is numerically nontrivial.
Let us first note the following generalization of \cite[Lemma 2.1]{Langer2024}:

\begin{Lemma}\label{short-exact-seq}
Let $$0\to \cE_1\to \cE \to \cE_2$$	be an exact sequence of reflexive sheaves on $X$ and assume that the cokernel of the last map is supported in codimension $\ge 2$. Then we have
	$$\int _X \ch _2 (\cE)L_1...L_{n-2} \le \int _X \ch _2 (\cE _1) L_1...L_{n-2}  +\int _X \ch _2 (\cE _2)L_1...L_{n-2} .$$
	Moreover, if the last map in the above sequence is surjective and the sequence is locally split in codimension $2$ then in the above inequality we have equality.
\end{Lemma}

\begin{proof} If $k$ has positive characteristic then the assertion follows from \cite[Lemma 2.1]{Langer2024}. So we can assume that  $k$ has characteristic $0$.
In case $X$ is a surface, the above inequality was proven in \cite[Lemma 1.2]{Langer2026WP}.
So the inequality follows from Theorem  \ref{high:properties} as in the proof of \cite[Lemma 2.1]{Langer2024}. For the last assertion, restrict to a geometric generic complete-intersection surface. Then the sequence becomes locally split everywhere, so Theorem \ref{surface:chern-ring}  gives additivity.
\end{proof}

\subsection{Bogomolov's inequality on normal projective varieties}

In this subsection and the next we assume that $k$ has characteristic $0$.
Let  $N_L (X)$ denote the quotient of the group of algebraic equivalence classes of Weil divisors on $X$  
by the radical of the intersection pairing
$$\langle D_1,D_2 \rangle _L:=D_1.D_2.L_2...L_{n-1}.$$
We define an open cone $K^+_L$ in $N_L (X)_{\RR}$ by 
$$K^+_L:=\{ D\in N_L (X)_{\RR}: \langle D, D \rangle _L>0 \hbox{ and }
\langle D, L_1 \rangle _L\ge 0\} .$$
By \cite[Lemma 1.20]{Langer2024},
$D\in K_L^{+}$ if and only if $\langle D, D' \rangle _L>0$ for all 
$D'\in {\overline{K_L^+}}\backslash \{ 0\} .$

Let us write $d$ for the number $L_1^2L_2...L_{n-1}$ and let us fix a coherent reflexive $\cO_X$-module $\cE$ of rank $r>0$. For a normal integral divisor $H\subset X$, we write $L_H$ for the
collection $(L_2|_H,\dots,L_{n-1}|_H)$.

 Note that Lemma \ref{short-exact-seq} implies that \cite[Lemma 2.2]{Langer2024} holds in characteristic zero.
So as in \cite{Langer2024} one can follow the proofs of \cite[Theorems 3.1, 3.2, 3.3 and 3.4]{Langer2004} to get the following results:

\begin{Theorem}
	Assume that $L_1$ is very ample and
	the restriction of $\cE$ to a very general divisor $H\in |L_1|$ is not slope
	$L_H$-semistable. 
	Let $r _i$ and $\mu_i$ denote  ranks and  slopes (with respect to $L_H$)
	of the Harder--Narasimhan filtration  of $\cE|_H$. Then
	$$\sum _{i<j} r_{i} r_{j}(\mu _{i} -\mu _{j})^2\le
	d \cdot \int_X\Delta (\cE)L_2\dots L_{n-1}+2r^2(\mu_{\max}(\cE) -\mu (\cE)) (\mu (\cE) -\mu_{\min} (\cE)) .
	$$ 
\end{Theorem}

In case all $L_i$ are ample  one can replace ``very general'' with ``general'' (see \cite{Langer2004-Add}).

\begin{Theorem}\label{Bog1}
	If $\cE$ is slope
	$L$-semistable  then 
	$ \int_X \Delta (\cE)L_2\dots L_{n-1}\ge 0.$
\end{Theorem}

\begin{Theorem}
	If $\int_X  \Delta (\cE) L_2\dots L_{n-1}<0$ 
	then there exists a rank $1\le r'<r$ saturated reflexive subsheaf $\cE '\subset \cE$ such that
	$\left( \frac{c_1(\cE')}{r'}- \frac{c_1(\cE)}{r}\right)$ lies in $ K^+_L.$ 
\end{Theorem}

The only difference in proofs with respect to \cite{Langer2004} is that we need Bogomolov's inequality on normal surfaces. This can be obtained by following classical Bogomolov's proof (see \cite{Bogomolov1979} or \cite[Theorem 3.4.1]{HL}). In fact, Theorem \ref{Bog1} follows also from \cite[Theorem 3.2]{Langer2004}
and Theorem \ref{main4} as noted in Corollary \ref{Bogomolov-inequality}.

As in \cite{Langer2024} and \cite{Langer2004} we get the following version of Bogomolov's inequality:

\begin{Theorem}
	If $\cE$ is a coherent reflexive $\cO_X$-module of rank $r$ then 
	$$d\cdot
	\int_X	\Delta (\cE)L_2\dots L_{n-1}+r^2(\mu_{\max} (\cE) -\mu (\cE)) (\mu (\cE) -\mu_{\min} (\cE) ) \ge 0.$$
\end{Theorem}

\subsection{Strong restriction theorem and boundedness}

Let us note that \cite[Lemma 2.3]{Langer2024} holds also in characteristic zero:

\begin{Lemma}
Let	$H\in|L_1|$ be a normal integral divisor and let $i:H\hookrightarrow X$ denote the closed embedding. 
Let $T$ be a torsion free coherent $\cO_H$-module of rank $\tau$ and let
$$0\to\cG\to\cE\to i_*T\to0$$
be a short exact sequence of coherent $\cO_X$-modules, where
$\cE$ is reflexive of rank $r$. Then $\cG$ is reflexive and
	\begin{align*}
		\int_X\Delta(\cG)L_2...L_{n-1}
		={}&\int_X\Delta(\cE)L_2...L_{n-1}
		-\tau(r-\tau)L_1^2L_2...L_{n-1}\\
		&+2\bigl(rc_1(T)-\tau c_1(i^*\cE)\bigr)
		.i^*L_2...i^*L_{n-1}.
	\end{align*}
\end{Lemma}

\begin{proof}
As in the proof of \cite[Lemma 2.3]{Langer2024} we can assume that $L_2,...,L_{n-1}$ are very ample. By Theorem~\ref{high:properties}, after extending the base field to an uncountable field, we can restrict to a very general complete intersection surface in $|L_2|\cap...\cap |L_{n-1}|.$	Such a surface is normal, the restrictions of $\cE$ and $\cG$ are reflexive, its intersection with  $H$ is a smooth curve and the restricted sequence remains exact with restriction of $T$ being torsion free. So we can reduce the assertion to the surface case.
	
In this case $H$ is a smooth Cartier curve contained in $X_{\mathrm{reg}}$, so $\cE$ and $\cG$ are isomorphic near every singular point of $X$, and their local Riemann--Roch correction terms coincide.
	Subtracting the Riemann--Roch formulas  gives
	\begin{align*}
		\chi(H,T)
		=\chi(X,\cE)-\chi(X,\cG)
		=\int_X\ch_2(\cE)-\int_X\ch_2(\cG)
		-\frac{\tau}{2}H.K_X.
	\end{align*}
	On the other hand, by the Riemann--Roch formula on $H$ we have
	$$\chi(H,T)=\deg_H T+\tau\chi(H,\cO_H)
	=\deg_H T-\frac{\tau}{2}(H^2+H.K_X).$$
	Consequently,
	$$\int_X\ch_2(\cE)-\int_X\ch_2(\cG)
	=\deg_H T-\frac{\tau}{2}H^2.$$
This implies that
	\begin{align*}
		\int_X\Delta(\cG)-\int_X\Delta(\cE)
		&=(c_1(\cE)-\tau H)^2-c_1(\cE)^2
		+2r\left(\deg_H T-\frac{\tau}{2}H^2\right)\\
		&=-\tau(r-\tau)H^2
		+2\bigl(r\deg_H T-\tau \deg_H(i^*\cE)\bigr),
	\end{align*}
which proves the assertion.
\end{proof}

The above lemma shows that the proof of \cite[Theorem 3.9]{Langer2024} works in characteristic zero and we have the following restriction theorem and its corollary:

\begin{Theorem}\label{strong-restriction-stable} 
	Let  $\cE$ be a coherent reflexive $\cO_X$-module of rank $r\ge 2$.  
Assume that $d>0$ and set 
	$$m_0=\frac{r-1}{ r}\int_X\Delta (\cE)L_2\dots L_{n-1} +\frac{1}{dr(r-1)}.$$
Assume that $\cE$ is slope $L$-stable and let $H\in |L_1^{\otimes m}|$  be an integral normal divisor.  
If $m>\left\lfloor m_0\right\rfloor$
then $\cE|_H$ is  slope $L_H$-stable.	
\end{Theorem}

\begin{Corollary}
	Let  $\cE$ be a coherent reflexive $\cO_X$-module of rank $r\ge 2$. Assume that $\cE$ is slope $L$-semistable 	and let $H\in |L_1^{\otimes m}|$  be an irreducible normal divisor. 
	Let $0= \cE_0\subset \cE_1\subset ...\subset \cE_s=\cE$ be a
	Jordan--H\"older filtration of $\cE$ and let us assume that all $(\cE_i/\cE_{i-1})|_H$ are torsion free. 
	Assume that $d>0$ and define $m_0$ as in Theorem \ref{strong-restriction-stable}. If 	$m>\left\lfloor m_0\right\rfloor$
	then $\cE|_H$ is slope $L_H$-semistable.
\end{Corollary}

This shows that the proofs of \cite[Theorem 3.14 and Corollaries 3.16, 3.17]{Langer2024} work also in characteristic zero. For simplicity we state just an analogue of \cite[Corollary 3.16]{Langer2024}.

\begin{Corollary} 
Assume that $L_1,...,L_{n-1}$ are ample and	
let us fix some positive integer $r$, integer $\ch_1$ and some real
	numbers $\ch _2$ and $\mu _{\max}$. 
	Let $\cB$ be the set of reflexive coherent $\cO_X$-modules $\cE$ such that 
	\begin{enumerate}
		\item $\cE$ has  rank $r$, 
		\item  $\int _X \ch_1 (\cE) L_1...L_{n-1}=\ch_1$,
		\item $\int _X \ch_2 (\cE) L_1...\widehat{L_i}...L_{n-1} \ge \ch_2$ for $i=1,...,n-1$,
		\item $\mu _{\max ,L } (\cE)\le \mu _{\max}$.
	\end{enumerate}
	Then the set $\cB$ is bounded.
\end{Corollary}

\subsection{A Miyaoka-type inequality}

Let $X$ be a normal projective variety of dimension $n\geq2$
over an algebraically closed field $k$.
For a collection $L=(L_1,...,L_{n-1})$ of nef line bundles (or just nef \(\mathbb Q\)-Cartier \(\mathbb Q\)-divisors) on $X$ we set
$$
L_{\min, L}(\cE)=
\begin{cases}
	\mu_{\min,L}(\cE),&\text{if }\operatorname{char}k=0,\\[4pt]
	\displaystyle\lim_{e\to\infty}
	\frac{\mu_{\min,L}(F_X^{[e]}\cE)}{p^e},
	&\text{if }\operatorname{char}k=p>0,
\end{cases}
$$
where $F_X^{[e]}\cE=((F_X^e)^*\cE)^{**}$.
We say that a
coherent reflexive sheaf $\cE$ is \emph{generically
	$L$-semipositive} if $L_{\min,L}(\cE)\geq0$.
We say that $\cE$ is \emph{generically semipositive} if it is
generically $L$-semipositive for every collection
$L$ of ample divisors.

The following result extends Miyaoka's theorem \cite[Theorem 6.1]{Miyaoka1985}
(see also \cite[Theorems~9.2 and~9.3]{LaGenericPos} for the positive characteristic versions) to reflexive sheaves on normal varieties.

\begin{Theorem}\label{Miyaoka-generically-nef}
	Let $\cE$ be a coherent reflexive sheaf on $X$, and assume
	that $L_1=c_1(\cE)$ is a nef $\QQ$-Cartier divisor.
	Fix nef line bundles $L_2,\ldots,L_{n-1}$.
	Assume that one of the following conditions holds:
	\begin{enumerate}
		\item The $1$-cycle $L_1\cdots L_{n-1}$ is numerically
		non-trivial and $\cE$ is generically
		$(L_1,\ldots,L_{n-1})$-semipositive.
		\item The $1$-cycle $L_1\cdots L_{n-1}$ is numerically
		trivial and $\cE$ is generically
		$(A,L_2,\ldots,L_{n-1})$-semipositive for some ample
		Cartier divisor $A$.
	\end{enumerate}
	Then
	$$
	\int_Xc_2(\cE)L_2\cdots L_{n-1}\geq0.
	$$
In particular, if \(\mathcal E\) is generically semipositive and \(c_1(\mathcal E)\) is nef and \(\mathbb Q\)-Cartier, then the inequality holds for every collection of nef line bundles \(L_2,\ldots,L_{n-1}\).
\end{Theorem}

\begin{proof}
	The numerical argument of \cite[Lemma~9.1]{LaGenericPos}
	applies to filtrations of reflexive sheaves with torsion free
	quotients, after replacing the quotients by their reflexive
	hulls. Indeed, Lemma~\ref{short-exact-seq} supplies the required
	inequality for $\ch_2$, and the Hodge index theorem applies
	to Weil-divisor intersections.
	
	In case~(1), the proof of \cite[Theorem~9.2]{LaGenericPos}
	therefore applies, using Corollary \ref{Bogomolov-inequality} and the ordinary Harder--Narasimhan
	filtration in characteristic zero and the strong Harder--Narasimhan filtration of a reflexive Frobenius	pullback in positive characteristic. In positive characteristic we also use Frobenius functoriality of the second Chern class.
	
	Consider case~(2).
Since $\mu_{(A,L_2,\ldots,L_{n-1})}(\cE)=0$, generic semipositivity implies that $\cE$ is  strongly $(A,L_2,\ldots,L_{n-1})$-semistable. Hence  if $A.L_2\cdots L_{n-1}$ is numerically
non-trivial, then Corollary \ref{Bogomolov-inequality} gives
	$$
	2r\int_Xc_2(\cE)L_2\cdots L_{n-1}
	\geq (r-1)L_1^2.L_2\cdots L_{n-1}=0,
	$$
where $r=\rk\cE$.	
	If $A.L_2\cdots L_{n-1}$ is numerically trivial, then
	$A^2.L_2\cdots L_{n-1}=0$. Nefness implies that 
	$L_2\cdots L_{n-1}$ has zero intersection with every
	codimension-two cycle. Singular Riemann--Roch therefore gives
	$$
	\chi\left(X,c_1(L_2)\cdots c_1(L_{n-1})\cdot\cF\right)=0
	$$
	for every coherent sheaf $\cF$.
	The asymptotic description of the second Chern class, together
	with $L_1^2.L_2\cdots L_{n-1}=0$, gives the asserted vanishing.
	
	For the last assertion, first suppose that \(L_2,\ldots,L_{n-1}\) are
	ample. For any ample divisor $A$, generic semipositivity
	with respect to $(L_1+\epsilon A,L_2,\ldots,L_{n-1})$ for
	all rational $\epsilon>0$ implies generic semipositivity
	with respect to $(L_1,\ldots,L_{n-1})$, by taking limits
	of the degrees of torsion free quotients, also after every
	reflexive Frobenius pullback in positive characteristic.
	Thus either~(1) or~(2) applies. The assertion for nef \(L_2,\ldots,L_{n-1}\)
	follows by continuity.
\end{proof}

\begin{Corollary}\label{Miyaoka-cotangent}
	Let $X$ be a normal projective variety of dimension $n\geq2$
	over an algebraically closed field $k$.
	Assume that $K_X$ is nef and $\QQ$-Cartier and that
	$\Omega_X^{[1]}$ is generically semipositive.
	Then for any nef line bundles $L_1,\ldots,L_{n-2}$ on $X$
	we have
	$$
	\int_Xc_2(\Omega_X^{[1]})L_1\cdots L_{n-2}\geq0.
	$$
	In characteristic zero, the generic semipositivity assumption
	can be replaced by the assumption that $X$ is non-uniruled.
\end{Corollary}

\begin{proof}
	The first assertion follows from
	Theorem~\ref{Miyaoka-generically-nef}.	
	For the second assertion we can assume that  $k=\mathbb C$.
	Let $\Omega_X^{[1]}\to\cQ$ be a nonzero torsion free quotient
	and let $f:Y\to X$ be a projective resolution which is an
	isomorphism over $X_{\mathrm{reg}}$.
	Saturating the extension of the kernel gives a torsion free
	quotient $\widetilde{\cQ}$ of $\Omega_Y^1$ with
	$f_*c_1(\widetilde{\cQ})=c_1(\cQ)$.
	Since $Y$ is non-uniruled,
	\cite[Theorem~0.1]{CampanaPeternell2011} implies that
	$\det\widetilde{\cQ}$ is pseudoeffective.
	Consequently, for any ample line bundles $L_1,\ldots,L_{n-1}$
	on $X$,
	$$
	c_1(\cQ).L_1\cdots L_{n-1}
	=
	c_1(\widetilde{\cQ}).f^*L_1\cdots f^*L_{n-1}\geq0.
	$$
	This proves generic semipositivity of $\Omega_X^{[1]}$.
\end{proof}

\begin{Remark}
If $X$ is a smooth complex projective variety with $K_X$ nef (or just pseudo-effective) then it is non-uniruled by
\cite{BDPP2013}. However, this fails for general normal complex projective varieties (e.g., for a cone over a smooth projective quintic curve in $\PP^2$).
\end{Remark}

\subsection{Big reflexive sheaves and hyperbolicity}

Let $\cE$ be a reflexive sheaf of rank $r$ on a normal proper surface $X$.
We say that $\cE$ is \emph{big} if
\[
\limsup_{m\longrightarrow\infty}
\frac{h^0(X,\Sym^{[m]}\cE)}{m^{r+1}}>0.
\]

The following theorem generalizes  \cite[Corollary 10.11]{Bogomolov1979}.

\begin{Theorem}\label{application:segre-bigness}
	Let $X$ be a normal projective surface over an algebraically closed
	field $k$. Let $\cE$ be a reflexive	sheaf of rank $r\geq1$ on $X$. Assume that $\cE$ is generically
	$H$-semipositive for some ample $H$. Then for every line bundle $A$ on $X$ we have
	\[
	h^0\bigl(X,\Sym^{[m]}\cE\otimes A^{-1}\bigr)
	\geq \frac{\int_Xs_2(\cE)}{(r+1)!}\,m^{r+1}-O(m^r).
	\]
	In particular, if $\int_Xs_2(\cE)>0$ then $\cE$ is big.
\end{Theorem}

\begin{proof}
	By Corollary \ref{sym:global-rr} we have	
	\[
	\chi(X,\Sym^{[m]}\cE\otimes A^{-1})= \chi(X,\Sym^{[m]}\cE)+O(m^r)
	=\frac{\int_X s_2(\cE)}{(r+1)!}\,m^{r+1}+O(m^r).
	\]	
	So it is sufficient to prove that $h^2(X,\Sym^{[m]}\cE\otimes A^{-1})= O(m^r)$. But by Serre's duality we have
	$h^2(X,\Sym^{[m]}\cE\otimes A^{-1})=h^0(X,(\Sym^{[m]}\cE)^*\hat \otimes A(K_X)).$
By the Ramanan--Ramanathan theorem (see \cite[Appendix]{Langer-Duke}) we have
\begin{align*}
\mu _{\max, H} ((\Sym^{[m]}\cE)^*\hat\otimes A(K_X))\le& L_{\max, H} ((\Sym^{[m]}\cE)^*\hat\otimes A(K_X))\\
&=-mL_{\min, H} (\cE)+(K_X+A)H\le (K_X+A)H.
\end{align*}
So the required inequality follows from the Le Potier--Simpson theorem (see \cite[Theorem 3.3.1]{HL})
and its positive characteristic version (see \cite[Corollary 3.4]{Langer-Duke}).
\end{proof}

As a corollary, we obtain the following result, which in the smooth case
was proven by Bogomolov in \cite{Bogomolov1977b}.

\begin{Corollary}\label{application:big-reflexive-cotangent}	
	Let $X$ be a normal projective surface over an algebraically closed
	field of characteristic zero. If $X$ is non-uniruled and \(\int_Xs_2(\Omega_X^{[1]})>0	\)
	then $\Omega_X^{[1]}$ is big and there are only finitely many complete integral curves
	\(
	C\subset X_{\mathrm{reg}},
	\)
	whose normalizations have genus at most one.
\end{Corollary}

\begin{proof}
	Since $X$ is non-uniruled, $\Omega_X^{[1]}$ is generically semipositive
	(see the proof of Corollary \ref{Miyaoka-cotangent}). Hence it is big
	by Theorem \ref{application:segre-bigness}. In particular, there exist
	an ample line bundle $A$, an integer $m>0$, and a nonzero section
	\[
	s\in H^0\!\left(
	X_{\mathrm{reg}},\operatorname{Sym}^{m}\Omega_{X_{\mathrm{reg}}}^{1}\otimes A^{-1}
	\right).
	\]
Then finiteness can be proven using Bogomolov's arguments as in	\cite[\S 5 and Theorem~2]{Bogomolov1977b} (see also	\cite[Proposition~4.3 and Theorem~4.5]{Deschamps1979}).
\end{proof}

The following example shows that the assumptions  of Corollary \ref{application:big-reflexive-cotangent}	 do not imply that a smooth
projective resolution of $X$ has a big cotangent bundle.

\begin{Example}
Let $B$ be a complex abelian surface and let $H$ be a very ample divisor on $B$. Fix an integer $k\geq2$ and	choose smooth curves $C\in|kH|$ and $D\in|(k+1)H|$ meeting transversely. 
Let $\pi:Y\to B$ be the blow up of $B$ at the points of $C\cap D$.
Let $E$ be the exceptional divisor and let $C'$	be the strict transform of $C$. Then
\(	L=\pi^*((k+1)H)-E\) is globally	generated and Stein's factorization of the morphism defined by the linear system $|L|$ gives a projective birational morphism
\( f:Y\to X	\)
	onto a normal projective surface, contracting exactly $C'$. Let $A$ be an ample line bundle on $X$
such that $L=f^*A$.	
Since $K_B=0$, we have	\(	K_Y= E=kL-(k+1)C'. \)
	The contracted point is singular, since $C'$ has positive genus.
	Restricting to $Y\setminus C'\simeq X_{\mathrm{reg}}$ and
	taking reflexive extensions shows that
	$\omega_X\simeq A^{\otimes k}$. In particular, $K_X$ is
	Cartier and ample, and the discrepancy of $C'$ is $-(k+1)$.
So we have
	\begin{align*}
		\int_Xs_2(\Omega_X^{[1]})
		&=\int_Ys_2(\Omega_Y^1)
		-\frac34c_1(f,\Omega_Y^1)^2+\frac14\Delta(f,\Omega_Y^1)\ge \frac{k(k+1)(3k-5)}4H^2>0.
	\end{align*}
Since $X$  is birational to $B$, it is non-uniruled and hence Corollary~\ref{application:big-reflexive-cotangent}
implies that $\Omega_X^{[1]}$ is big.
	
	On the other hand, the symmetric plurigenera of smooth
	projective varieties are birational invariants and
	$\Omega_B^1\simeq\cO_B^{\oplus2}$, so
	$h^0(Y,\Sym^m\Omega_Y^1)=m+1$
	and the cotangent bundle of $Y$, or of any smooth projective
	resolution of $X$, is not big. 
\end{Example}

\medskip

The following example coming from \cite[Theorem 5.1]{Garbagnati-Penegini} shows that in general, the assumptions of 
Corollary \ref{application:big-reflexive-cotangent}	 do not imply that there are only finitely many integral curves \(C\subset X,\) whose normalizations have genus at most one.

\begin{Example}
Consider the smooth complex projective genus-two curve $C$
with an affine part given by \(y^2=x^5+1\)	and the automorphism \(\gamma(x,y)=(\zeta_5x,y)\) of $C$, where $\zeta_5$ is a primitive fifth root of unity. Let us consider the surface
	\(
	X:=(C\times C)/\langle\gamma\times\gamma^2\rangle.
	\)
	It is normal projective with cyclic quotient singularities. The quotient map \(q:C\times C\to X\) is quasi-étale of degree \(5\), and \(K_X\) is ample.
Since \(q^{[*]}\Omega_X^{[1]}\simeq\Omega_{C\times C}^1, \)
	we have
	\[
	K_X^2=\frac85> 
	\int_Xc_2(\Omega_X^{[1]})=\frac45.
	\]	
 By \cite[Theorem 5.1]{Garbagnati-Penegini},
the smooth minimal model of $X$ is a K3 surface.
In particular, $X$ is non-uniruled. The arguments in \cite[Section 7.2]{Garbagnati-Penegini} show that $X$
admits a covering family of integral curves of geometric genus one.
More precisely,  let $x_i,y_i$ denote the pullbacks of the affine coordinates $x,y$
under the projections $C\times C\to C$, for $i=1,2$.
The rational functions $t=x_2/x_1^2$, $y_1$ and $y_2$ are invariant
under $\gamma\times\gamma^2$, and they give the presentation
\[
\CC(X)\simeq
\operatorname{Frac}\left(
\CC[t,y_1,y_2]/
\bigl(y_2^2-t^5(y_1^2-1)^2-1\bigr)
\right).
\]
The polynomial $t^5(y_1^2-1)^2+1\in\CC(t)[y_1]$ has degree four
and four distinct roots over $\overline{\CC(t)}$.
Thus the smooth projective model of the generic fibre of
the rational map $t:X\dashrightarrow\mathbb P^1$ is a double
cover of $\mathbb P^1_{\CC(t)}$ branched at four geometric points.
By the Riemann--Hurwitz formula, it has genus one.
Moreover, $(y_1,y_2)=(1,1)$ gives a $\CC(t)$-rational point.
The closures in $X$ of the general fibres therefore form
a covering family of integral curves whose normalizations
are elliptic curves.

Note however that \(X_{\mathrm{reg}}\) is Brody hyperbolic as an entire curve in \(X_{\mathrm{reg}}\) lifts through the finite étale cover \(q^{-1}(X_{\mathrm{reg}})\to X_{\mathrm{reg}}\), and both projections of its lift to \(C\) are constant. Thus every elliptic curve in the covering family necessarily meets the singular locus of $X$.
\end{Example}

\appendix
\section{Non-existence of a good theory of higher Chern classes}

The second Chern classes (or Chern characters) constructed in this paper are compatible with
finite pullbacks and the asymptotic Riemann--Roch formula. In higher degrees, the
expected leading term in asymptotic Riemann--Roch is already
incompatible with the usual duality formula, even without imposing
any functoriality under pullbacks. The obstruction occurs for
rank-one reflexive sheaves on a projective threefold.

For a rank-one reflexive sheaf $\cL$ we write
$\cL^{[m]}=(\cL^{\otimes m})^{**}$.

\begin{Proposition}\label{app:nonexistence}
	Let $k$ be an algebraically closed field and let $n\geq3$.
	There is no assignment of a real number $\gamma_X(\cL)$ to every
	normal projective $n$-dimensional variety $X$ over $k$ and every
	isomorphism class of rank-one reflexive sheaves $\cL$ on $X$
	satisfying the following properties:
	\begin{enumerate}
		\item For every rank-one reflexive sheaf $\cL$,
		$\gamma_X(\cL^*)=(-1)^n\gamma_X(\cL).$
		\item For every rank-one reflexive sheaf $\cL$,
		$\chi(X,\cL^{[m]})
		=\gamma_X(\cL^{[m]})+O(m^{n-1}).$
	\end{enumerate}
\end{Proposition}

\begin{proof}
	Since $(\cL^{[m]})^*\simeq(\cL^*)^{[m]}$, the two properties imply
	$$
	\chi(X,\cL^{[m]})
	-(-1)^n\chi(X,(\cL^*)^{[m]})
	=O(m^{n-1}).
	$$
Consider the projective quadric cone
$X=\{x_0x_3-x_1x_2=0\}\subset\PP^4_k$
and  set $D=\{x_0=x_1=0\}\subset X$.
One can easily check that for every integer $m\geq0$ we have
\begin{equation*}
	\chi(X,\cO_X(mD))=m+1,
	\qquad
	\chi(X,\cO_X(-mD))=\binom{m+1}{3}-m+1.
\end{equation*}	
So for $n=3$ and $\cL=\cO_X(D)$ we have
	$$
	\chi(X,\cL^{[m]})+\chi(X,(\cL^*)^{[m]})
	=\binom{m+1}{3}+2,
	$$
	contradicting the required $O(m^2)$ estimate.
For $n>3$, one can obtain a similar contradiction considering
$X\times\PP^{n-3}$ and $\cL\boxtimes\cO_{\PP^{n-3}}(1).$
\end{proof}

The intended invariant would be
$\gamma_X(\cL)=\int_X\ch_n(\cL)$.
Thus, in every degree $n\geq3$, the usual duality formula
and the expected leading term in asymptotic Riemann--Roch
cannot hold simultaneously for all rank-one reflexive sheaves.
 The proof also works with an error
$o(m^n)$ in property~(2).

\begin{Remark}
	A related failure of duality for pullbacks of Weil divisors was
	observed by de Fernex and Hacon
	\cite[Remark 2.11 and Example 2.12]{dFH09}.
	For a birational morphism $f:Y\to X$ of normal varieties, their
	natural pullback is characterized by
	$$
	\cO_Y(-f^\natural D)
	=\bigl(\cO_X(-D)\cdot\cO_Y\bigr)^{**}.
	$$
	They show that $f^\natural(-D)\geq-f^\natural D$ and that
	their asymptotic pullback can satisfy
	$f^*(-D)\ne-f^*D$.
	Their example exhibits a failure of pullback constructions
	to respect duality. Proposition~\ref{app:nonexistence}
	concerns the compatibility of duality with asymptotic
	Euler characteristics.
\end{Remark}

\section*{AI Declaration}

During the development and preparation of this article, the author used ChatGPT (OpenAI; GPT-6 Astra) extensively as an interactive tool to explore proof strategies, discuss mathematical arguments,  and assist with drafting and revision. 
The author independently checked the mathematical arguments and cited sources, rewrote proofs and the text, and takes full responsibility for the results and the final manuscript.

\section*{Acknowledgements}

The author was partially supported by the National Science Centre, Poland, contract number 2025/59/B/ST1/02168.

\bibliographystyle{amsalpha-addendum}
\bibliography{References}
\end{document}